%% file: spinsl2-ii.tex
\documentclass[12pt]{amsart}

\usepackage{amsmath}%
\usepackage{amsfonts}%
\usepackage{amssymb}%
\usepackage{graphicx}
\usepackage[all]{xy}
\usepackage{array}
\usepackage{url}

\usepackage{tikz}
\usetikzlibrary{decorations.shapes}
\usetikzlibrary{decorations.text}
\usetikzlibrary{decorations.pathreplacing}
\usetikzlibrary{decorations.markings}
\usetikzlibrary{shapes}
\usetikzlibrary{shapes.symbols}
\input{tikzstyledefs}
\tikzstyle{every picture}=[semithick,baseline=0pt,heightone,label distance=-2mm]

\newtheorem{theorem}{Theorem}
\theoremstyle{plain}

\newtheorem{corollary}[theorem]{Corollary}

\newtheorem{definition}[theorem]{Definition}
\newtheorem{example}{Example}

\newtheorem{lemma}[theorem]{Lemma}

\newtheorem*{notation*}{Notation}

\newtheorem{proposition}[theorem]{Proposition}

\numberwithin{equation}{section}

\numberwithin{theorem}{section}
\numberwithin{equation}{section}

\def\C{\mathbb{C}}
\def\na{\underline{a}}

\def\tr{\mathrm{tr}}

\def\sixj#1#2#3#4#5#6{\begin{bmatrix}#1 & #2 & #6\\#4 & #3 & #5\end{bmatrix}}
\def\sixjt#1#2#3#4#5#6{\bigl[\begin{smallmatrix}#1 & #2 & #6\\#4 & #3 & #5\end{smallmatrix}\bigr]}

\def\edges#1#2#3{\mathfrak{e}_{#1}(#2,#3)}

\def\mfe{\mathfrak{e}}
\def\mfb{\mathfrak{b}}
\def\mff{\mathfrak{f}}
\def\mfs{\mathfrak{s}}

\def\Fus#1#2#3#4#5{\mathfrak{F}_{#1}\left(\begin{smallmatrix}#3 & #5\\ #2 & #4\end{smallmatrix}\right)}
\def\nFus#1#2#3#4#5{\hat{\mathfrak{F}}_{#1}\left(\begin{smallmatrix}#3 & #5\\ #2 & #4\end{smallmatrix}\right)}

\newcommand{\Tr}[1]{\mathrm{tr}(#1)}
\newcommand{\xb}{\mathbf{X}}
\newcommand{\yb}{\mathbf{Y}}
\newcommand{\zb}{\mathbf{Z}}

\newcommand{\ub}{\mathbf{U}}
\newcommand{\vb}{\mathbf{V}}
\newcommand{\wb}{\mathbf{W}}

\newcommand{\aq}{/\!\!/}
\newcommand{\X}{\mathfrak{X}}
\newcommand{\Y}{\mathfrak{Y}}
\newcommand{\R}{\mathfrak{R}}
\newcommand{\hm}{\mathrm{Hom}}
\newcommand{\SL}{\mathrm{SL}(2,\C)}
\newcommand{\SLm}[1]{\mathrm{SL}(#1,\C)}

\newcommand{\F}{\mathtt{F}}

\newcommand{\xt}{\mathtt{x}}

\newcommand{\id}{\mathbf{I}}

\newcommand{\End}{\mathrm{End}}
\newcommand{\iadm}[1]{\lceil #1 \rfloor}
\newcommand{\vi}{\vec{i}}
\newcommand{\vj}{\vec{j}}
\newcommand{\vk}{\vec{k}}

\newcommand{\tmx}[2]{\left(\begin{matrix}#1\\#2\end{matrix}\right)}                      
\newcommand{\tmxt}[2]{\bigl(\begin{smallmatrix}#1\\#2\end{smallmatrix}\bigr)} 
\newcommand{\Adm}[3]{#3\in\lceil #1,#2\rfloor}
\newcommand{\N}{\mathbb{N}}
\newcommand{\ch}[3]{\chxx_{\sss #2,#3}^{\sss #1}}
\newcommand{\chh}[3]{\chxx^{\sss #2,#3}_{\sss #1}}
\def\sss{\scriptscriptstyle}
\newcommand{\chxx}{\raisebox{2pt}{$\chi$}}
\newcommand{\cspan}{\C\,}

\newcommand{\bs}[3]{{#1}^*_{#2}\otimes {#1}_{#3}}
\def\vprod#1{V^{\otimes #1}}

\newcommand{\imx}[1]{\tmx{#1_{11} & #1_{12}}{#1_{21} & #1_{22}}}

\begin{document}
\title[Computing Central Functions]{Computing $\mathrm{SL}(2,\C)$ Central Functions with Spin Networks}
\author{Sean Lawton and Elisha Peterson}
\address[Sean Lawton]{University of Texas, Pan American, 1201 West University Drive, Edinburg, TX 78539}
\email[]{lawtonsd@utpa.edu}
\address[Elisha Peterson]{Johns Hopkins University Applied Physics Lab, Laurel, MD 20723}
\email[]{elisha.peterson@jhuapl.edu}
\date{\today}

\begin{abstract}
Let $G=\SLm{2}$ and $\F_r$ be a rank $r$ free group.  Given an {\it admissible} weight $\vec{\lambda}$ in $\mathbb{N}^{3r-3}$, there exists a class function defined on $\hm(\F_r,G)$ called a {\it central function}.   We show that these functions admit a combinatorial description in terms of graphs called {\it trace diagrams}.  We then describe two algorithms (implemented in {\it Mathematica}) to compute these functions.
\end{abstract}
\maketitle


\section{Introduction}\label{intro-section}\label{s:introduction}
\input{intro.tex}

\section{Spin Networks and Trace Diagrams}\label{s:spinnets}
\input{spinnetsintro.tex}

\section{Rank $r$ Central Functions}\label{s:centralfunctions}
\input{rankrcentral.tex}

\section{Trace Diagram Recurrences}\label{s:recurrences}
\input{tracediagramreview.tex}
\input{simplerecurrence.tex}

\input{nonsimplerecurrence.tex}

\section{Rank Three Central Functions}\label{s:rankthree}
\input{rank3.tex}

\section*{Acknowledgments}
We would like to thank Bill Goldman for introducing us to this topic.  A very special thanks is due to Suhyoung Choi and to Korea Advanced Institute of Science and Technology (KAIST) for generously inviting us to speak for their seventh Geometric Topology Fair.  This paper comes out of our time at KAIST. We also thank the referee for several valuable comments, including particularly helpful suggestions for clarifying the statements of Theorem \ref{t:simplerecurrence} and Corollary \ref{c:simplerecurrence}.

\appendix
\section{Background}\label{s:background}
\input{appendix}

\bibliographystyle{amsalpha}
\bibliography{bib}

\end{document}

%% file: tikzstyledefs.tex
\tikzstyle heightone=[scale=.7,shift={(0,-.3)}]
\tikzstyle heightones=[scale=.8,xscale=.35,shift={(0,.1)}]
\tikzstyle heightoneonehalf=[scale=.9,shift={(0,-.2)}]
\tikzstyle heighttwo=[scale=.9,shift={(0,-.4)}]
\tikzstyle heighttwos=[scale=.5,xscale=.6,shift={(0,-.1)}]
\tikzstyle heightthree=[scale=.6,shift={(0,-.9)}]
\tikzstyle heightthrees=[scale=.4,xscale=.7,shift={(0,-.2)}]

\tikzstyle vertex=[circle,draw,fill=black,inner sep=1pt]
\tikzstyle ciliation=[circle,draw=none,fill=red,inner sep=1pt,semitransparent]
\tikzstyle ciliatednode=[vertex,pin={[pin distance=1mm,pin edge={semitransparent,red},ciliation]#1:{}}]

\tikzstyle matrix=[black,thick,circle,draw=blue!50,top color=blue!20,bottom color=black!10,scale=.8,inner sep=1pt]
\tikzstyle small matrix=[matrix,scale=.7]
\tikzstyle vector=[black,thick,rectangle,draw=gray!50!yellow,top color=yellow!30,bottom color=black!10,scale=.8,inner sep=2pt]
\tikzstyle small vector=[vector,scale=.8]
\tikzstyle plain vector=[rectangle,draw=none,fill=white,scale=.7]

\tikzstyle basiclabel=[draw=none,fill=none,shape=rectangle,inner sep=2pt,scale=.8]
\tikzstyle leftlabel=[basiclabel,anchor=east]
\tikzstyle rightlabel=[basiclabel,anchor=west]
\tikzstyle bottomlabel=[basiclabel,anchor=north]
\tikzstyle toplabel=[basiclabel,anchor=south]

\tikzstyle trivalent=[very thick]

\tikzstyle arrowstyle=[blue,semitransparent,scale=2]

\tikzstyle directed=[postaction={decorate,decoration={markings,
    mark=at position .65 with {\arrow[arrowstyle]{stealth}}}}]
\tikzstyle reverse directed=[postaction={decorate,decoration={markings,
    mark=at position .65 with {\arrowreversed[arrowstyle]{stealth};}}}]

\tikzstyle with matrix=[postaction={decorate,decoration={markings,
    mark=at position .5 with {\node[matrix]{#1};}}}]
\tikzstyle with small matrix=[postaction={decorate,decoration={markings,
    mark=at position .5 with {\node[small matrix]{#1};}}}]

\tikzstyle directed matrix=[postaction={decorate,decoration={markings,
    mark=at position .9 with {\arrow[arrowstyle]{stealth}},
    mark=at position .35 with {\node[matrix]{#1};}}}]
\tikzstyle directed small matrix=[postaction={decorate,decoration={markings,
    mark=at position .9 with {\arrow[arrowstyle]{stealth}},
    mark=at position .35 with {\node[small matrix]{#1};}}}]
\tikzstyle reverse directed matrix=[postaction={decorate,decoration={markings,
    mark=at position .4 with {\arrowreversed[arrowstyle]{stealth};},
    mark=at position .65 with {\node[matrix]{#1};}}}]
\tikzstyle reverse directed small matrix=[postaction={decorate,decoration={markings,
    mark=at position .4 with {\arrowreversed[arrowstyle]{stealth};},
    mark=at position .65 with {\node[small matrix]{#1};}}}]

\tikzstyle dotdotdot=[decorate,decoration={markings,
    mark=at position .3 with{\node{.};},
    mark=at position .5 with {\node{.};},
    mark=at position .7 with {\node{.};}}]

\tikzstyle wavyup=[out=90,in=-90]
\tikzstyle wavydown=[out=-90,in=90]

%% file: intro.tex
Let $\F_r=\langle \xt_1,...,\xt_r\rangle$ be a rank $r$ free group and let $G=\SL$.  The representations $\R_r=\hm(\F_r,G)$ constitute an affine variety isomorphic to $G^{\times r}$.  $G$ acts on $\R_r$ by conjugation.

As an affine variety $\R_r$ is associated to a reduced $\C$-algebra $\C[\R_r]$, and as $G$-modules we have the following decomposition:
$$\C[\R_r]^G \approx \sum_{\vec{i}\in\N^r} \sum_{\substack{ \vec{j},\vec{k} \in \iadm{\vec{i}} \\ |\vec{k}|=|\vec{j}| }} \End
\left(V_{\left(|\vec{i}|-2|\vec{j}|\right)}\right)^{G}$$
where $V_k=\mathrm{Sym}^k(\C^2)$, $\vec{i}=(i_1,i_2,...,i_r)\in \N^r$, $|\vec{i}|=i_1+\cdots + i_r$, and $\vec{j}=(j_1,...,j_{r-1})\in \N^{r-1}$ satisfies $\vec{j}\in \iadm{\vec{i}}$ if and only if 
$$0\leq j_l \leq \mathrm{min}(i_1+\cdots + i_l -2(j_1 +\cdots + j_{l-1}),i_{l+1})$$ for all $l=1,2,...,r-1$.

This generalizes the complexification of the Peter-Weyl Theorem ${\displaystyle \C[G]^G\approx \sum_{k\in \N}\End(V_k)^G}$.

\begin{definition}
Given the above isomorphism, for each triple $\vec{i},\vec{j},\vec{k}$ such that $\vec{i}\in\N^r$, $\vec{j},\vec{k}\in \iadm{\vec{i}}$, and $|\vec{j}|=|\vec{k}|$,  there exists a class function $\chh\vi\vj\vk\in\mathrm{End}(V_{\left(|\vec{i}|-2|\vec{j}|\right)})^{G}\subset \C[\R_r]^G.$  We refer to the functions $\chh\vi\vj\vk$ as rank $r$ \emph{central functions}.
\end{definition}

Given any two central functions $\chi$ and $\chi'$, their product $\chi\chi'$ is again a polynomial function in $\C[\R_r]^G$, and thus $\chi\chi'= \sum \lambda_k \chi_k$ for coefficients $\lambda_k$ and other central functions $\chi_k$.  This ring structure is very far from the natural ring structure of $\C[\R_r]^G\subset \C[\R_r]$, and is quite mysterious.  It seems, thus far, to be best understood in terms of special graphs called {\it trace diagrams}.

In \cite{LP} we described rank 1 and rank 2 central functions and determined their ring structure by exploring the calculus of trace diagrams.  It is the purpose of this paper to show that this association generalizes and to describe algorithms to compute these functions explicitly in the case of rank 3 central functions.

\begin{figure}[htb]
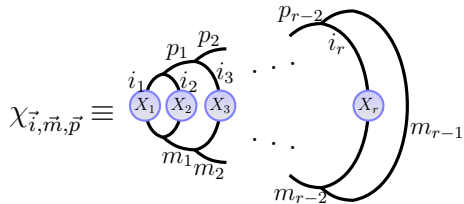

    $$
    \chi_{\vi,\vec m,\vec p} \equiv
    \tikz[scale=1.2]{
        \draw[trivalent]
            (0,0)to[bend left=80]node[small matrix]{$X_1$}(0,1)node[leftlabel,pos=.8]{$i_1$}
            (0,0)to[bend right=80]node[small matrix]{$X_2$}(0,1)node[rightlabel,pos=.8]{$i_2$}
            (0,0)to[bend right=20](.5,-.2)node[bottomlabel,pos=.5]{$m_1$}
            to[bend right=80]node[small matrix]{$X_3$}(.5,1.2)node[rightlabel,pos=.8]{$i_3$}
            to[bend right=20](0,1)node[toplabel,pos=.5]{$p_1$}
            (.5,-.2)to[bend right=20](1,-.4)node[bottomlabel,pos=.5]{$m_2$}
            (1,1.4)to[bend right=20](.5,1.2)node[toplabel,pos=.5]{$p_2$};
        \draw[draw=none](1.25,0)--(2.25,-.2)node[pos=.2]{.}node[pos=.5]{.}node[pos=.8]{.};
        \draw[draw=none](1.25,1)--(2.25,1.2)node[pos=.2]{.}node[pos=.5]{.}node[pos=.8]{.};
        \draw[trivalent,shift={(.5,0)}]
            (1.5,-.6)to[bend right=20](2,-.8)node[bottomlabel,pos=.4]{$m_{r-2}$}
            to[bend right=80]node[small matrix]{$X_r$}(2,1.8)node[leftlabel,pos=.8]{$i_r$}
            to[bend right=20](1.5,1.6)node[toplabel,pos=.6]{$p_{r-2}$}
            (2,-.8)to[bend right=20](2.5,-1)
            to[bend right=80](2.5,2)node[rightlabel,pos=.4]{$m_{r-1}$}
            to[bend right=20](2,1.8);
    }$$
    \caption{Diagrammatic Form of Central Functions}\label{f:centralfunction}
    \end{figure}

The diagramatic form, shown in Figure \ref{f:centralfunction}, can be deduced to exist from work of Baez \cite{Ba96} and later Sikora \cite{Si2}. One of the most important points of this paper is that central functions can be understood entirely in terms of combinatorics given by a graphical calculus.  This culminates with Theorem \ref{simplelooptheorem}.

The paper is organized as follows.  In the next section, we introduce spin networks and trace diagrams.  Then we go through the construction of rank $r$ central functions and work out a few examples in section 3.  In section 4 we further the study of the trace diagram calculus in preparation for the algorithm to compute them.  In section 5 we give two algorithms to compute rank 3 central functions that we have successfully implemented in {\it Mathematica}. In the appendix, we review some relevant invariant theory and representation theory.

%% file: spinnetsintro.tex
The tools used in this paper to explore central functions are \emph{trace diagrams}, which are a slight generalization of \emph{spin networks}, which will be formally defined in Definitions \ref{def:spinnet} and \ref{def:tracediagram}.

Informally, a spin network is a graph that can also be interpreted as a function. Depending upon the type of diagram, the function may be between tensor products of $V=\C^2$, or between tensor products of $\SL$ representations, which are symmetric powers of $V$. In order for a spin network's function to be well-defined, the inputs and outputs of the function must be specified, and the graph must be given a small amount of extra structure. Once this is done, algebraic rules become local diagrammatic rules that are much easier to work with.

Trace diagrams generalize spin networks by allowing matrices to act on specific strands of the diagram. When the diagram has no inputs or outputs, it becomes a function $\SL\times\cdots\times\SL\to\C$ that is invariant under simultaneous conjugation, hence an element of $\C[\R_r]^G$ for some $r$.

\begin{figure}[htb]
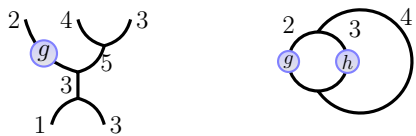

    $$\tikz[trivalent,shift={(0,.5)}]{
        \draw(-.5,-1)node[leftlabel]{$1$}to[bend left](0,-.5);
        \draw(.5,-1)node[rightlabel]{$3$}to[bend right](0,-.5);
        \draw(0,-.5)to[wavyup]node[leftlabel]{$3$}(0,0)to[bend left]node[matrix]{$g$}(-1,1)node[leftlabel]{$2$};
        \draw(0,0)to[bend right]node[rightlabel]{$5$}(.5,.5)to[bend right](1,1)node[rightlabel]{$3$};
        \draw(.5,.5)to[bend left](0,1)node[leftlabel]{$4$};
    }
    \qquad\qquad
    \tikz[trivalent,every node/.style={basiclabel},scale=1.4]{
        \draw(0,.5)circle(.4)(0,.1)arc(-145:145:.7);
        \node[matrix]at(-.4,.5){$g$};
        \node[matrix]at(.4,.5){$h$};
        \node at(-.4,1){$2$};\node at(.5,.95){$3$};\node at(1.2,1.1){$4$};
    }
    $$
\caption{Example of a spin network representing a function from $V_1\otimes V_3$ to $V_2\otimes V_3\otimes V_4$, and a trace diagram representing a function $\SL\times\SL\to\C$.}\label{f:spinnet-tracediagram-example}
\end{figure}

Diagrammatic notations such as spin networks were developed by physicists in the middle of the 20th century as a graphical description of quantized angular momentum. The theory has been developed extensively since that time, most notably by Cvitanovic \cite{Cv08} and Stedman \cite{ste90}.  Trace diagram calculations are very similar to skein module calculations in knot theory. Indeed, the Kauffman bracket skein module can be thought of as a quantization of the diagrams in this paper \cite{Kau91}.

This section follows the diagrammatic conventions of \cite{LP,Pet06}, rather than the more common conventions of \cite{CFS95,Kau91,Ma99}. The most obvious difference is \eqref{eq:capcup}, which means that several diagrams will differ in sign from these sources. Proofs are omitted, as they do not contribute to the main point of the paper. We refer the reader to \cite{LP,Pet06} for these. Similar proofs for alternate diagrammatic conventions are in \cite{CFS95,Kau91}.

\subsection{Diagrammatic Description of $V_n$}\label{ss:diagrammatic-symmetrization}

We represent vectors $v\in V=\C^2$ by $\tikz{\draw(0,.2)node[vector]{$v$}to(0,.8);}$, and dual vectors $v^*\in V^*$ by $\tikz{\draw(0,.2)to(0,.8)node[vector]{$v$};}$. Strands without vector labels are unknown vectors. Tensor products are represented by drawing diagrams adjacent to each other:
    $v_1\otimes v_2\equiv\tikz{\draw(0,.2)node[small vector]{$v_1$}to(0,.8)(.5,.2)node[small vector]{$v_2$}to(.5,.8);}$.
The symmetrization operation is indicated by
    \begin{equation}\label{eq:symmetrizer-diagram}
    \tikz[xscale=.7]{
        \foreach\xa in{1,1.3,2.7,3}{\draw(\xa,0)to(\xa,1);}
        \draw[draw=black,fill=gray!10](.8,.35)rectangle(3.2,.75);\node[basiclabel]at(2,.55){$n$};
        \node[basiclabel]at(2,.1){$\cdots$};
    }
    =\frac{1}{n!}\sum_{\sigma\in\Sigma_n}
    \tikz[xscale=.7]{
        \foreach\xa in{1,1.3,2.7,3}{\draw(\xa,0)to(\xa,1);}
        \draw[draw=black,fill=gray!10](.8,.35)rectangle(3.2,.75);\node[basiclabel]at(2,.55){$\sigma$};
        \node[basiclabel]at(2,.1){$\cdots$};
    },
    \end{equation}
which is an \emph{idempotent} since
    $$\tikz[xscale=.7]{
        \foreach\xa in{1,1.3,2.7,3}{\draw(\xa,0)to(\xa,1);}
        \draw[draw=black,fill=gray!10](.8,.35)rectangle(3.2,.75);\node[basiclabel]at(2,.55){$n$};
        \node[basiclabel]at(2,.1){$\cdots$};
    }\circ\tikz[xscale=.7]{
        \foreach\xa in{1,1.3,2.7,3}{\draw(\xa,0)to(\xa,1);}
        \draw[draw=black,fill=gray!10](.8,.35)rectangle(3.2,.75);\node[basiclabel]at(2,.55){$n$};
        \node[basiclabel]at(2,.1){$\cdots$};
    }=\tikz[xscale=.7]{
        \foreach\xa in{1,1.3,2.7,3}{\draw(\xa,0)to(\xa,1);}
        \draw[draw=black,fill=gray!10](.8,.35)rectangle(3.2,.75);\node[basiclabel]at(2,.55){$n$};
        \node[basiclabel]at(2,.1){$\cdots$};
    }.$$

The basis elements for $V_n=\mathrm{Sym}^n(V)$ and $V_n^*=\mathrm{Sym}^n(V^*)$ described in section \ref{ss:reptheory} take the form
    \begin{equation}\label{eq:symmetric-basis}
    {\sf n}_{n-k}=
    \tikz[xscale=.7]{
        \foreach\xa/\xb in{1/$e_1$,3/$e_1$,4/$e_2$,6/$e_2$}{
            \draw(\xa,0)node[small vector]{\xb}to(\xa,1);
        }
        \draw[draw=black,fill=gray!10](.8,.45)rectangle(6.2,.8);\node[basiclabel]at(3.5,.625){$n$};
        \node[basiclabel]at(2,0){$\cdots$};\node[basiclabel]at(5,0){$\cdots$};
        \draw[decorate,decoration=brace](3.2,-.35)to node[auto,basiclabel]{$n-k$}(.8,-.35);
        \draw[decorate,decoration=brace](6.2,-.35)to node[auto,basiclabel]{$k$}(3.8,-.35);
    }
    \qquad
    {\sf n}^*_{n-k}=
    \tikz[xscale=.7,yscale=-1,shift={(0,-.5)}]{
        \foreach\xa/\xb in{1/$e_1$,3/$e_1$,4/$e_2$,6/$e_2$}{
            \draw(\xa,0)node[small vector]{\xb}to(\xa,1);
        }
        \draw[draw=black,fill=gray!10](.8,.45)rectangle(6.2,.8);\node[basiclabel]at(3.5,.625){$n$};
        \node[basiclabel]at(2,0){$\cdots$};\node[basiclabel]at(5,0){$\cdots$};
        \draw[decorate,decoration=brace](.8,-.35)to node[auto,basiclabel]{$n-k$}(3.2,-.35);
        \draw[decorate,decoration=brace](3.8,-.35)to node[auto,basiclabel]{$k$}(6.2,-.35);
    }.
    \end{equation}

It is also customary to represent elements of $V_n$ and $V_n^*$ by thicker strands labeled by $n$, so that the basis elements just described may also be written
    \begin{equation}\label{eq:thick-strands}
    {\sf n}_i = \tikz[trivalent]{\draw(0,0)node[small vector]{${\sf n}_i$}to(0,1)node[rightlabel]{$n$};}
    \qquad
    {\sf n}^*_i = \tikz[trivalent]{\draw(0,0)node[rightlabel]{$n$}to(0,1)node[small vector]{${\sf n}_i$};}
    \end{equation}

The pairing between $V_n$ and $V_n^*$ is given by
    $${\sf n}^*_{n-k}({\sf n}_{n-l})
    =\tikz[trivalent]{\draw(0,-.2)node[small vector]{${\sf n}_{n-l}$}to(0,1.2)
        node[rightlabel,pos=.6,scale=.8]{$n$}node[small vector]{${\sf n}_{n-k}$};}
    =\raisebox{2pt}{$\delta_{kl}$}\!\Big/\!\raisebox{-2pt}{$\tbinom{n}{k}$}.$$

The $\SL$-action on $V_n$ is represented by a matrix $g\in\SL$ drawn on the strand:
    $$g\cdot {\sf n}_{n-k}
    = \tikz[trivalent]{\draw(0,-.2)node[small vector]{${\sf n}_{n-k}$}
        to(0,.6)node[matrix]{$g$}to(0,1.2)node[rightlabel]{$n$};}$$

\subsection{Spin Networks}
The formal definitions of spin networks and trace diagrams follow.

\begin{definition}\label{def:spinnet}
A \emph{spin network} is a trivalent graph drawn in the plane with edges labeled by representations of $\mathrm{SL}(2,\mathbb{C})$, such that the labels adjoining each vertex form an admissible triple. The graph is drawn with all free ends at the top or the bottom, and all local extrema are in general position relative to the orientation of the diagram.
\end{definition}

\begin{definition}\label{def:tracediagram}
A \emph{trace diagram} is a spin network with one or more edges marked by elements of the matrix group $\mathrm{SL}(2,\mathbb{C})$, with no marking occurring at a local extrema.
\end{definition}

Every trace diagram corresponds to a unique function. The domain of the function is a tensor product of irreducible $\SL$-representations, given by the labeling of the input edges, while the co-domain is given by the labeling of the output edges.

The function may be computed by decomposing the diagram into horizontal slices. The diagram's function is the composition of the functions of the individual slices. (One can also decompose the diagram into vertical slices, which are conjoined as tensor products.)

We now describe the lexicon for reading off a spin network's function.  Let $\{e_1, e_2\}$ be the standard basis for $\C^2$.  The single-strand \emph{cup} and \emph{cap} diagrams are
    \begin{equation}\label{eq:capcup}
    \tikz{\draw(0,.7)arc(180:360:.5);} = e_1\otimes e_2-e_2\otimes e_1
    \qquad\text{and}\qquad
    \tikz{\draw(0,.2)arc(180:0:.5);} : v\otimes w \mapsto \det[v\: w],
    \end{equation}
respectively. Building upon this definition, the symmetrized cup and cap maps are
    \begin{equation}\label{eq:trivalent-cap}
    \tikz[trivalent]{\draw(0,.7)arc(180:360:.5)node[rightlabel]{$n$};}\equiv
    \tikz{\draw(-.1,1)--(-.1,.8)arc(180:360:.3)--(.5,1);\draw(-.7,1)--(-.7,.8)arc(180:360:.9)--(1.1,1);
        \draw[draw=black,fill=gray!10](.05,.7)rectangle(.35,-.2);\node[basiclabel,scale=.9]at(.2,.25){$n$};
        \node[basiclabel,scale=.7]at(-.375,.8){$\cdots$};\node[basiclabel,scale=.7]at(.825,.8){$\cdots$};}
    \qquad\text{and}\qquad
    \tikz[trivalent]{\draw(0,.2)arc(180:0:.5)node[rightlabel]{$n$};}\equiv
    \tikz[yscale=-1,shift={(0,-.8)}]{\draw(-.1,1)--(-.1,.8)arc(180:360:.3)--(.5,1);\draw(-.7,1)--(-.7,.8)arc(180:360:.9)--(1.1,1);
        \draw[draw=black,fill=gray!10](.05,.7)rectangle(.35,-.2);\node[basiclabel,scale=.9]at(.2,.25){$n$};
        \node[basiclabel,scale=.7]at(-.375,.8){$\cdots$};\node[basiclabel,scale=.7]at(.825,.8){$\cdots$};}
    .
    \end{equation}

From \eqref{eq:capcup} and the action of $\SL$ on $V=\C^2$, one can deduce 
    \begin{equation}\label{eq:trace-diagram}
        \tikz{\draw(0,.5)circle(.5);\draw(.5,.5)node[matrix]{$g$};} = \tr(g).
    \end{equation}
This is why the diagrams are referred to as \emph{trace diagrams}.

Using Schur's Lemma, the Clebsch-Gordan decomposition \eqref{eq:clebsch-gordan} implies that for each admissible triple $\{a,b,c\}$, as specified in Definition \ref{def:admissible}, there are unique (up to a nonzero multiple) injections $V_c\hookrightarrow V_a\otimes V_b$ and projections $V_a\otimes V_b\twoheadrightarrow V_c$. These are sometimes called \emph{intertwining operators}. The injection has the following diagrammatic form:
    \begin{equation}\label{eq:trivalent-vertex}
    \tikz[trivalent,shift={(0,.5)}]{
        \foreach\xa/\xb in{150/a,30/b,-90/c}{\draw[rotate=\xa](0,0)--(0:.8)(0:1)node[basiclabel]{$\xb$};}}
    =
    \tikz[scale=.6,shift={(0,1)}]{
        \foreach\xa in{0,120,240}{
            \draw[rotate=\xa](28:2.2)to[bend right](-88:2.2)(0:2.5)to[bend right](-60:2.5);}
        \foreach\xa/\xb in{30/b,150/a,270/c}{
            \draw[rotate=\xa,draw=black,fill=gray!10](1.2,1.5)rectangle(1.8,-1.5);
            \draw[rotate=\xa](0:1.52)node[basiclabel]{$\xb$};
            \draw[rotate=\xa](60:1.25)node[basiclabel,rotate=\xa,rotate=60]{$\cdots$};
        }
    }
    \end{equation}
The projection $V_a\otimes V_b\twoheadrightarrow V_c$ is obtained by reflection through a horizontal line.

The following notation will be used extensively in later sections.
\begin{notation*}
    Given an admissible triple $\{a,b,c\}$, define
        $$ \mathfrak{e}_a(b,c) = \frac{b+c-a}{2}, \qquad \mathfrak{e}(a,b,c) = \frac{a+b+c}{2}.$$
\end{notation*}
In the figure $\mathfrak{e}_a(b,c)$ represents the (unique) number of strands connecting the $b$ and $c$ edges in the expansion. In these terms, the admissibility condition in Definition \ref{def:admissible} can be reformulated as
    $$\mathfrak{e}_a(b,c) \geq 0, \quad \mathfrak{e}_b(a,c) \geq 0, \quad \mathfrak{e}_c(a,b) \geq 0.$$

\subsection{Spin Network Relations}

Spin networks satisfy certain \emph{skein relations} that can be leveraged to reason about the underlying functions. The \emph{spin network skein module} is the space of formal sums of spin networks with coefficients in $\mathbb{C}$, modulo certain relations arising from the representation theoretic ``building blocks'' \eqref{eq:symmetrizer-diagram}-\eqref{eq:trivalent-vertex}.

The relations that follow can be proven directly from spin network definitions. Proofs are contained in \cite{LP} and earlier works.

A general spin network skein relation has the form $\sum_m \alpha_m \mathsf{s}_m = 0$, where $\alpha_m\in\mathbb{C}$ and each $\mathsf{s}_m$ represents a diagram. Such relations are well-behaved under reflections:

\begin{proposition}[Spin Network Reflection, Proposition 3.6 in \cite{Pet06}]\label{p:spinnetreflection}
    Given a spin network $\mathsf{s}$, denote by $\mathsf{s}^\updownarrow$ the spin network obtained by vertical reflection of $\mathsf{s}$, and by $\mathsf{s}^\leftrightarrow$ the spin network obtained by horizontal reflection of $\mathsf{s}$. Then
        $$
        \sum_m \alpha_m \mathsf{s}_m = 0
            \quad \Leftrightarrow \quad \sum_m \alpha_m \mathsf{s}^\updownarrow_m = 0
            \quad \Leftrightarrow \quad \sum_m \alpha_m \mathsf{s}^\leftrightarrow_m = 0.
        $$
\end{proposition}

Other kinds of topological moves may introduce signs. This is the case for local extrema introduced into a strand, crossings of edges adjacent to a vertex, and reorientation of vertices, as indicated in Figure \ref{f:signs}.
\begin{figure}[htp]
    $$
    \tikz{\draw(0,0)to[wavyup](.1,1);}
    =
    -\tikz{\draw(.1,0)to[out=90,in=-110](0,.5)to[out=70,in=-110,looseness=2](.4,.5)to[out=70,in=-90,looseness=1](.3,1);}
    \qquad
    \tikz[xscale=.6]{\draw(0,.2)arc(180:0:.5);}
    =
    -\tikz[shift={(0,.4)}]{\coordinate(vxa)at(0,.55)edge[bend right=45](-.3,.2)edge[bend left=45](.3,.2);
        \draw(-.3,-.4)to[wavyup](.3,.2)(.3,-.4)to[wavyup](-.3,.2);}
    $$
\caption{Diagrammatic moves on single strands in a spin network that introduce signs.} \label{f:signs}
\end{figure}

Define $\mfs_a(b,c) \equiv (-1)^{\mfe_a(b,c)}$.

\begin{proposition}[Spin Network Sign Changes, Proposition 3.22 in \cite{LP}]\label{p:spinnetsigns}
    \begin{align}
    \tikz[trivalent]{
        \draw(.1,0)to[out=90,in=-110](0,.5)to[out=70,in=-110,looseness=2](.4,.5)to[out=70,in=-90,looseness=1](.3,1)node[rightlabel]{$n$};}
    \label{eq:kink}
    & = (-1)^n \tikz[trivalent]{\draw(0,0)to[wavyup](.1,1)node[rightlabel]{$n$};}
    \\
    \tikz[trivalent]{
        \coordinate(vxa)at(0,.55)
            edge[bend right]node[leftlabel,pos=1]{$a$}(-.4,0)edge[bend left]node[rightlabel,pos=1]{$b$}(.4,0)
            edge[]node[rightlabel,pos=1]{$c$}(0,1);}
    \label{eq:crossedvertex}
    & = \mfs_c(a,b)
    \tikz[trivalent,shift={(0,.2)}]{
        \coordinate(vxa)at(0,.55)
            edge[bend right](-.3,.2)edge[bend left](.3,.2)edge[]node[rightlabel,pos=1]{$c$}(0,1);
        \draw(-.3,-.4)node[leftlabel]{$a$}to[wavyup](.3,.2)(.3,-.4)node[rightlabel]{$b$}to[wavyup](-.3,.2);}
    = \mfs_a(b,c)
    \tikz[trivalent]{
        \coordinate(vxa)at(0,.45)
            edge[]node[leftlabel,pos=1]{$a$}(0,0)
            edge[bend left](-.2,.7)edge[bend right](.2,.65);
        \draw(-.2,.7)to[wavyup](-.1,1.1)node[leftlabel]{$c$}
            (.45,0)node[rightlabel]{$b$}to[wavyup](.55,.6)to[bend right=90,looseness=2](.2,.65);}
    \end{align}
\end{proposition}
A more general version of this result follows. See Figure \ref{f:signs} for the meaning of `kink' and `crossed extrema'.
\begin{proposition}\label{p:spinnetsignstrong}
Let $\mathsf{T}_1$ and $\mathsf{T}_2$ be topologically equivalent spin networks (where the topological equivalence respects both the position of free ends and the labels on edges). Then
$$\mathsf{T}_1=(-1)^{k+x}\mathsf{T}_2,$$
where $k$ is the number of ``kinks'' and $x$ the number of ``crossed extrema'' in the diagrams obtained by expanding edges as in \eqref{eq:trivalent-vertex}.
\end{proposition}

The following two propositions describe the most basic diagrams, as well as how to join and separate strands in spin networks.

\begin{proposition}[Closed Spin Networks, Proposition 3.19 in \cite{LP}]\label{p:deltatheta}
    \begin{align*}
    \tikz[trivalent]{\draw(0,.5)circle(.5);\node[basiclabel]at(.5,1){$c$};}
    & \equiv \Delta(c) = c+1 = \mathsf{dim}(V_c); \\
    \tikz[trivalent,every node/.style={basiclabel}]{
        \draw(0,.5)circle(.4)(0,.1)arc(-145:145:.7);
        \node at(-.5,.85){$a$};\node at(.5,.85){$b$};\node at(1.2,1.1){$c$};}
    & \equiv \Theta(a,b,c) = \frac{\edges abc! \edges bac! \edges cab! (\mfe(a,b,c)+1)!}{a!b!c!}.
    \end{align*}
\end{proposition}

\begin{proposition}[Bubble/Fusion Identities, Propositions 3.20,3.21 in \cite{LP}]\label{p:bubblefusion}
    \begin{align}
    \tikz[trivalent,every node/.style={basiclabel}]{
        \draw(0,.5)circle(.3);\node at(-.5,.7){$a$};\node at(.5,.7){$b$};
        \draw(-.05,-.2)node[anchor=west]{$c$}to[wavyup](0,.2)(0,.8)to[wavyup](.05,1.2)node[anchor=west]{$d$};
    } \label{eq:bubbleidentity}
    & = \mfb_c(a,b) \: \tikz[trivalent]{\draw(-.05,-.2)to[wavyup](.05,1.2)node[rightlabel]{$c$};} \! \delta_{c,d}
    \\
    \tikz[trivalent]{
        \draw(-.3,-.1)to[wavyup](-.2,.5)to[wavyup](-.3,1.1)node[leftlabel]{$a$}
            (.3,-.1)to[wavyup](.2,.5)to[wavyup](.3,1.1)node[rightlabel]{$b$};
    } \label{eq:fusionidentity}
    & = \sum_{c\in\iadm{a,b}} \mff_c(a,b)
    \tikz[trivalent]{
    \coordinate(vxa)at(0,.25)
        edge[bend right]node[leftlabel,pos=1]{$a$}(-.3,-.15)edge[bend left]node[rightlabel,pos=1]{$b$}(.3,-.15);
    \coordinate(vxb)at(0,.75)
        edge[]node[basiclabel,auto]{$c$}(vxa)
        edge[bend left]node[leftlabel,pos=1]{$a$}(-.3,1.15)edge[bend right]node[rightlabel,pos=1]{$b$}(.3,1.15);},
    \end{align}
    where $\mfb_c(a,b) = \frac{\Theta(a,b,c)}{\Delta(c)} = \frac{1}{\mff_c(a,b)}$.
\end{proposition}
We refer to $\mfb_c(a,b)$ and $\mff_c(a,b)$ as the \emph{bubble constant} and \emph{fusion constant}, respectively.

\subsection{$\SL$-Equivariance and Trace Diagrams}
The rules for manipulating trace diagrams are the same as those for manipulating spin networks, since all parts of a diagram, as mappings, are $\SL$-equivariant. In particular:
\begin{proposition}[Equivariance of Trace Diagrams, Propositions 3.8 and 3.11 in \cite{LP}]\label{s:diagram-equivarance}
    For $g\in\SL$,
    \begin{align*}
        \tikz{\draw(0,.7)arc(180:360:.5);} &= 
            \tikz{\draw(0,1)--(0,.7)arc(180:360:.5)--(1,1);
            \draw(0,.85)node[matrix]{$g$};\draw(1,.85)node[matrix]{$g$};}\\
        \tikz[xscale=.7]{
            \foreach\xa in{1,1.3,2.7,3}{\draw(\xa,0)to(\xa,1.3)node[small matrix,pos=.3]{$g$};}
            \draw[draw=black,fill=gray!10](.8,.7)rectangle(3.2,1.1);\node[basiclabel]at(2,.9){$n$};
            \node[basiclabel]at(2,.1){$\cdots$};
        }&=\tikz[xscale=.7]{
            \foreach\xa in{1,1.3,2.7,3}{\draw(\xa,0)to(\xa,1.3)node[small matrix,pos=.7]{$g$};}
            \draw[draw=black,fill=gray!10](.8,.2)rectangle(3.2,.6);\node[basiclabel]at(2,.4){$n$};
            \node[basiclabel]at(2,1.2){$\cdots$};
        }\\
        \tikz[trivalent]{\draw(0,.7)arc(180:360:.5)node[rightlabel]{$n$};} &=
            \tikz[trivalent]{\draw(0,1)--(0,.7)arc(180:360:.5)--(1,1);\draw(.95,.4)node[rightlabel]{$n$};
            \draw(0,.85)node[matrix]{$g$};\draw(1,.85)node[matrix]{$g$};}\\
        \tikz[trivalent,shift={(0,1)},yscale=-1]{
            \coordinate(vxa)at(0,.55)
            edge[bend right]node[leftlabel,pos=1]{$a$}(-.6,-.4)edge[bend left]node[rightlabel,pos=1]{$b$}(.6,-.4)
            edge[]node[rightlabel,pos=1]{$c$}(0,1.1);
            \draw(-.5,0)node[matrix]{$g$};\draw(.5,0)node[matrix]{$g$};
            } &=
        \tikz[trivalent,shift={(0,1.3)},yscale=-1]{
            \coordinate(vxa)at(0,.55)
            edge[bend right]node[leftlabel,pos=1]{$a$}(-.6,-.1)edge[bend left]node[rightlabel,pos=1]{$b$}(.6,-.1)
            edge[]node[rightlabel,pos=1]{$c$}(0,1.4);
            \draw(0,1)node[matrix]{$g$};
            }.
    \end{align*}
\end{proposition}
Note the consequences
    $$
        \tikz{\draw(0,.7)arc(180:360:.5)--(1,1);\draw(0,.85)node[matrix]{$g$};} =
            \tikz{\draw(0,1)--(0,.7)arc(180:360:.5)--(1,1);
            \draw(1,.85)node[small matrix,ellipse]{$g^{-1}$};}
        \qquad\text{and}\qquad
        \tikz[trivalent]{\draw(0,.7)arc(180:360:.5)--(1,1)node[rightlabel]{$n$};\draw(0,.85)node[matrix]{$g$};} =
            \tikz[trivalent]{\draw(0,1)--(0,.7)arc(180:360:.5)--(1,1);\draw(.95,.4)node[rightlabel]{$n$};
            \draw(1,.85)node[small matrix,ellipse]{$g^{-1}$};}.
    $$

These equivariance properties permit every relation of the previous section to be adapted to include matrices along the edges. In later sections, we will often display diagrammatic manipulations without the required matrices, appealing to equivariance to show that they also hold with the matrices.

%% file: rankrcentral.tex
\subsection{Decomposition Theorem}
Let $\C[\SL]$ be the coordinate ring of the variety $G=\SL$, and recall our convention that $V_k$ denotes the $k^{\rm th}$ symmetric tensor of $\C^2$, denoted by $\mathrm{Sym}^k(\C^2)$.

The following theorem is a consequence of the ``unitary trick'', the Peter-Weyl Theorem, and the fact that the set of
matrix coefficients of $\SL$ is exactly its coordinate ring. See \cite{LP} for a detailed proof.

\begin{theorem}[Decomposition]\label{decomposition}
There is an $\SL$-module isomorphism $$\C[\SL]\approx\sum_{k\in \N}V_k^*\otimes V_k\approx \sum_{k\in \N}\End(V_k).$$
\end{theorem}

The isomorphism is given by defining $$\Upsilon :\sum_{n\ge0} V_n^*\otimes V_n\longrightarrow \C [\SL]$$ by linear extension of
the mapping
$${\sf n}^*_{n-k}\otimes {\sf n}_{n-l}\mapsto {\sf n}^*_{n-k}(\xb\cdot {\sf n}_{n-l}),$$ where
$\xb=\imx{x}$ is a generic matrix.

In particular, letting $\{e_1, e_2\}$ be the standard basis for $\C^2$,
\begin{align}\label{eq:tensormatrixcontraction}
    {\sf n}^*_{n-k}(\xb\cdot {\sf n}_{n-l})
    = \tikz[trivalent,heighttwo]{
        \draw(0,0)node[small vector]{${\sf n}_{n-l}$}
            to node[small matrix]{$\mathbf{X}$}node[pos=.75,rightlabel]{$n$}(0,2)node[small vector]{${\sf n}_{n-k}$};}
    &=  {\sf n}^*_{n-k}\left((x_{11} e_1+x_{21} e_2)^{n-l}(x_{12} e_1 + x_{22} e_2)^l\right)\nonumber\\
    &=\sum_{\substack{i+j=k\\0 \le i \le n-l \\ 0 \le j \le l}}
        \tbinom{n}{k}^{-1}\tbinom{n-l}{i}\tbinom{l}{j}x_{11}^{n-l-i}x_{12}^{l-j}x_{21}^ix_{22}^j.
\end{align}

\subsection{Applying the decompostion.}

\begin{eqnarray*}
\C[\SL^{\times r}] & \approx& \C[\SL]^{\otimes r}\\
 & \approx& \bigotimes_{1\leq k\leq r}\left(\sum_{i_k\in \N}V^*_{i_k}\otimes V_{i_k}\right)\\
 & \approx& \sum_{(i_1,...,i_r)\in\N^r} V_{i_1}^*\otimes V_{i_1}\otimes \cdots \otimes V_{i_r}^* \otimes V_{i_r}\\
 & \approx& \sum_{(i_1,...,i_r)\in\N^r} V_{i_1}^*\otimes \cdots \otimes V^*_{i_r}\otimes V_{i_1}\otimes \cdots \otimes V_{i_r} \ .\\
\end{eqnarray*}

As is stated in \cite{Pet06}, the isomorphism above is determined by the following association:  $(v^*_{i_1}\otimes v^*_{i_2}\otimes \cdots \otimes v^*_{i_r} )\otimes (v_{j_1}\otimes v_{j_2}\otimes \cdots \otimes v_{j_r} )$ maps to the polynomial function
$$(\xb_1, \xb_2,...,\xb_r) \mapsto v^*_{i_1}(\xb_1\cdot v_{j_1})v^*_{i_2}(\xb_2\cdot v_{j_2})\cdots v^*_{i_r} (\xb_r\cdot v_{j_r} ).$$  We will call this ``tensorial contraction.''

Our principal interest is with the {\it invariant} polynomial functions that arise in this fashion.  To determine these polynomials we will need a notion of ``admissibility.''

We say $(\{V_{i_1},...,V_{i_r}\}, V_x)$ is an admissible pair if and only if $V_x$ occurs as a summand in the decomposition of $V_{i_1}\otimes\cdots \otimes V_{i_r}$ into irreducible $G$-modules.  In this case there exists a $G$-module $W$ for which $V_{i_1}\otimes\cdots \otimes V_{i_r}\approx V_x\oplus W$ (as $G$-modules). The existence of an injection corresponds to a way to connect a single $x$ strand to the $i_1,i_2,\ldots,i_r$ strands in an admissible way:
    $$\tikz[heighttwo]{
        \node[cloud,fill=red!20,draw=red!20!gray,inner sep=3pt](middle)at(0,1){??}
            edge[trivalent]node[rightlabel,pos=1]{$x$}(0,-.5)
            edge[trivalent,bend left=10]node[rightlabel,pos=1]{$i_1$}(-1.2,2.5)
            edge[trivalent,bend left=10]node[rightlabel,pos=1]{$i_2$}(-.5,2.5)
            edge[trivalent,draw=none]node[basiclabel]{$\cdots$}(.35,2.5)
            edge[trivalent,bend right=10]node[rightlabel,pos=1]{$i_r$}(1,2.5);
    }
    $$
Several injections are possible, but in this paper we focus on the \emph{left-associative} injection
    \begin{equation}\label{eq:leftassocdiagram}
    \tikz[heightthree,scale=1.5]{
        \draw[trivalent](0,0)node[rightlabel]{$x$}to(0,1)to[bend left]node[leftlabel]{$m_{r-2}$}(-.5,1.5)
            (0,1)to[bend right](1.5,3.5)node[rightlabel]{$i_r$}
            (-1,2)to[bend left]node[leftlabel]{$m_2$}(-1.5,2.5)to[bend left]node[leftlabel]{$m_1$}(-2,3)
                to[bend left](-2.3,3.5)node[leftlabel]{$i_1$}
            (-2,3)to[bend right](-1.75,3.5)node[rightlabel]{$i_2$}
            (-1.5,2.5)to[bend right](-1,3.5)node[rightlabel]{$i_3$};
        \draw[dotdotdot,bend right](-.7,3)to(.7,2.5);
    }
    \end{equation}
This diagram is only admissible if the triples at each vertex are admissible, meaning $m_1\in\iadm{i_1,i_2}$, and for $l>1$, $m_l\in\iadm{m_{l-1},i_{l+1}}$.

From Proposition \ref{clebshgordan}, $\{V_{i_1},V_{i_2}\}$ and $V_x$ are admissible if and only
if $x=i_1+i_2-2j_1$ for $0\leq j_1\leq \mathrm{min}(i_1,i_2)$.

Now consider $\{V_{i_1},V_{i_2},V_{i_3}\}$ and $V_x$.  Using the above example and Proposition \ref{clebshgordan} a second time we have $$V_x\hookrightarrow (V_{i_1}\otimes V_{i_2}) \otimes V_{i_3}\approx \sum V_{i_1+i_2-2j_1}\otimes V_{i_3}\approx \sum  V_{i_1+i_2+i_3-2(j_1+j_2)}$$ where $0\leq j_1\leq \mathrm{min}(i_1,i_2)$ and $0\leq j_2\leq \mathrm{min}(i_1+i_2-2j_1,i_3)$.  Therefore, $(\{V_{i_1},V_{i_2},V_{i_3}\},V_x)$ is an admissible pair whenever $x=i_1+i_2+i_3-2(j_1+j_2)$ and both inequalities $0\leq j_1\leq \mathrm{min}(i_1,i_2)$ and $0\leq j_2\leq \mathrm{min}(i_1+i_2-2j_1,i_3)$ are satisfied.

Generalizing these examples by iteratively using the Clebsch-Gordan formula to decompose $V_{i_1}\otimes\cdots\otimes V_{i_r}$ we come to the following notation and definition.

Let $\vec{i}=(i_1,i_2,...,i_r)\in \N^r$, and let $|\vec{i}|=i_1+\cdots + i_r$.
\begin{definition}
    We say that $\vec{j}=(j_1,...,j_{r-1})\in \N^{r-1}$ is $\vec{i}$-admissible $($and denote it by $\vec{j}\in \iadm{\vec{i}})$ if and only if for all $1\leq l \leq r-1$ we have
    $$0\leq j_l \leq \mathrm{min}(i_1+\cdots + i_l -2(j_1 +\cdots + j_{l-1}),i_{l+1}).$$
\end{definition}
Note that this is precisely the condition $m_i\in\iadm{i_1,i_2}, \: m_l\in\iadm{m_{l-1},i_{l+1}}$ given earlier with $m_l=i_1+\cdots+i_{l+1}-2(j_1+\cdots+j_n)$.

Recall $G=\SL$.  We then use Clebsch-Gordon iteratively with respect to Theorem \ref{decomposition} to conclude:
\begin{eqnarray*}
\C[G^{\times r}] & \approx& \C[G]^{\otimes r}\\
 & \approx& \sum_{(i_1,...,i_r)\in\N^r} V_{i_1}^*\otimes \cdots \otimes V^*_{i_r}\otimes V_{i_1}\otimes \cdots \otimes V_{i_r}\\
 & \approx& \sum_{\vec{i}\in\N^r} \left( \sum_{\vec{j}\in \iadm{\vec{i}}} V_{\left(|\vec{i}|-2|\vec{j}|\right)}^*\right)\otimes \left(\sum_{\vec{k}\in
\iadm{\vec{i}}} V_{\left(|\vec{i}|-2|\vec{k}|\right)}\right)\\
& \approx& \sum_{\vec{i}\in\N^r} \sum_{\vec{j},\vec{k} \in \iadm{\vec{i}}} V_{\left(|\vec{i}|-2|\vec{j}|\right)}^*\otimes V_{\left(|\vec{i}|-2|\vec{k}|\right)}\ .\\
\end{eqnarray*}
Since the above maps are $\SL$-equivariant,
\begin{equation*}
 \C[\R_r]^G=\C[G^{\times r}]^{G} \approx \sum_{\vec{i}\in\N^r} \sum_{\vec{j},\vec{k} \in \iadm{\vec{i}}} \left(V_{\left(|\vec{i}|-2|\vec{j}|\right)}^*\otimes
V_{\left(|\vec{i}|-2|\vec{k}|\right)}\right)^{G}.
\end{equation*}

By Schur's Lemma,
$$\mathrm{dim}_{\C}\left(V_{\left(|\vec{i}|-2|\vec{j}|\right)}^*\otimes V_{\left(|\vec{i}|-2|\vec{k}|\right)}\right)^{G}
= \left\{
\begin{array}{ll} 1 & \textrm{if $|\vec{k}|=|\vec{j}|$}\\ 0 & \textrm{if $|\vec{j}|\not=|\vec{k}|$}\end{array}\right.$$
Therefore,

$$\C[\R_r]^G=\C[G^{\times r}]^{G} \approx \sum_{\vec{i}\in\N^r} \sum_{\substack{ \vec{j},\vec{k} \in \iadm{\vec{i}} \\ |\vec{k}|=|\vec{j}| }} \End
\left(V_{\left(|\vec{i}|-2|\vec{j}|\right)}\right)^{G}.$$

\begin{definition}
Given the above isomorphism, for each triple $\vec{i},\vec{j},\vec{k}$ such that $\vec{i}\in\N^r$, $\vec{j},\vec{k}\in \iadm{\vec{i}}$, and $|\vec{j}|=|\vec{k}|$,  there exists a class function $\chh\vi\vj\vk\in\C[\R_r]^G$ which corresponds to a generating homothety $($unique up to scalar$)$ in $\mathrm{End}(V_{\left(|\vec{i}|-2|\vec{j}|\right)})^{G}$. We refer to the functions $\chh\vi\vj\vk$ as \emph{rank $r$ central functions}.
\end{definition}

Denote by $\cspan\chh\vi\vj\vk \subset \C[\R_r]^G$ the linear span over $\C$ of $\chh\vi\vj\vk$.

In these terms,
$$\C[\R_r]^G\approx \sum_{\vec{i}\in\N^r} \sum_{\substack{\vec{j},\vec{k} \in \iadm{\vec{i}} \\ |\vec{k}|=|\vec{j}|}} \cspan\chh\vi\vj\vk.$$

Thus, the central functions $\chh\vi\vj\vk$ form an additive basis for $\C[\R_r]^G$.  However, the multiplicative structure in terms of this basis is very complicated and not at all obvious.

We note that $\vi$ has $r$ entries, $\vk$ and $\vj$ have $r-1$.  So the index relation $|\vj|=|\vk|$ shows that each central function is in terms of exactly $3r-3$ indices, the Krull dimension
of $\C[\R_r]^G$.  Thus each central function corresponds to an admissible weight $\vec{\lambda}\in \mathbb{N}^{3r-3}$.

Let the Clebsch-Gordan injection be denoted by $$\iota^{\sss \vk}_{\sss \vi}:V_{\left(|\vec{i}|-2|\vec{k}|\right)}\hookrightarrow V_{i_1}\otimes\cdots \otimes V_{i_r}.$$
Also, let $\{{\sf c}^*_s\}$ be a basis for $V^*_{\left(|\vec{i}|-2|\vec{j}|\right)}$ and $\{{\sf d}_t\}$ is a basis for $V_{\left(|\vec{i}|-2|\vec{k}|\right)}$ (assuming $|\vj|=|\vk|$).

Define $$\mathcal{M}_{\vi}^{\vj,\vk}=\bigg(\iota^{\sss \vj}_{\sss \vi}({\sf c}_s^*)\Big((\xb_1,...,\xb_r)\cdot\iota^{\sss \vk}_{\sss \vi}({\sf
d}_t)\Big)\bigg)_{\!\!st}.$$  $\mathcal{M}_{\vi}^{\vj,\vk}$ is a $\left(|\vec{i}|-2|\vec{k}|+1\right) \times \left(|\vec{i}|-2|\vec{k}|+1\right)$ matrix with noted $s,t$ entries.

In these terms we can see that

$$\chh\vi\vj\vk(\xb_1,...,\xb_r)=\Tr{\mathcal{M}_{\vi}^{\vj,\vk}}.$$

Since these injections are given by iteratively using the injections from the rank 2 case (that is decomposing a product of tensors two at a time), our computation of these injections in \cite{LP} determine all such injections in general (up to a choice of associativity).

With this in mind, these functions take natural diagrammatic form. Beginning with \eqref{eq:leftassocdiagram} and its vertical reflection (providing the decomposition of the dual), tensorial contraction corresponds to gluing copies of the matrix variables $\xb_l$ in between the two diagrams. Taking the trace corresponds to adding a closing loop to the diagram. The resulting diagram is
    \begin{equation}\label{eq:rankrdiagram}
    \chh\vi\vj\vk(\xb_1,\ldots,\xb_r) \equiv
    \chi_{\vi,\vec m,\vec p} \equiv
    \tikz[scale=1.2]{
        \draw[trivalent]
            (0,0)to[bend left=80]node[small matrix]{$X_1$}(0,1)node[leftlabel,pos=.8]{$i_1$}
            (0,0)to[bend right=80]node[small matrix]{$X_2$}(0,1)node[rightlabel,pos=.8]{$i_2$}
            (0,0)to[bend right=20](.5,-.2)node[bottomlabel,pos=.5]{$m_1$}
            to[bend right=80]node[small matrix]{$X_3$}(.5,1.2)node[rightlabel,pos=.8]{$i_3$}
            to[bend right=20](0,1)node[toplabel,pos=.5]{$p_1$}
            (.5,-.2)to[bend right=20](1,-.4)node[bottomlabel,pos=.5]{$m_2$}
            (1,1.4)to[bend right=20](.5,1.2)node[toplabel,pos=.5]{$p_2$};
        \draw[draw=none](1.25,0)--(2.25,-.2)node[pos=.2]{.}node[pos=.5]{.}node[pos=.8]{.};
        \draw[draw=none](1.25,1)--(2.25,1.2)node[pos=.2]{.}node[pos=.5]{.}node[pos=.8]{.};
        \draw[trivalent,shift={(.5,0)}]
            (1.5,-.6)to[bend right=20](2,-.8)node[bottomlabel,pos=.4]{$m_{r-2}$}
            to[bend right=80]node[small matrix]{$X_r$}(2,1.8)node[leftlabel,pos=.8]{$i_r$}
            to[bend right=20](1.5,1.6)node[toplabel,pos=.6]{$p_{r-2}$}
            (2,-.8)to[bend right=20](2.5,-1)
            to[bend right=80](2.5,2)node[rightlabel,pos=.4]{$m_{r-1}$}
            to[bend right=20](2,1.8);
    }
    \end{equation}
where $m_l=i_1+\cdots+i_{l+1}-2(j_1+\cdots+j_l)$ and $p_l=i_1+\cdots+i_{l+1}-2(k_1+\cdots+k_l)$. The requirement $|\vj|=|\vk|$ becomes $m_{r-1}=p_{r-1}$.

\subsection{Example $r=1$}
The diagram is a single loop:
    $$
    \chi_c=
    \tikz[trivalent]{
        \draw(0,.5)circle(.5);
        \node[basiclabel]at(.5,1){$c$};
        \node[small matrix]at(-.5,.5){$\mathbf{X}$};
    }
    $$

The trivial representation $V_0$ gives $\chxx_0=1$. The standard representation $V_1$ has diagonal matrix coefficients $x_{11}$ and $x_{22}$, hence
$$\chxx_1= x_{11}+x_{22}=\Tr{\xb}.$$

The remaining functions may be computed directly, or by using the following product formula:
    \begin{equation}\label{eq:rank1product}
        \chxx_a \chxx_b = \sum_{c\in \lceil a,b\rfloor} \chxx_c
    \end{equation}
Explicitly, the particular case $b=1$ is (for $a\ge1$)
    \begin{equation}\label{eq:rank1recurrence}
        \chxx_a \chxx_1 = \chxx_{a+1} + \chxx_{a-1},
    \end{equation}
from which the recurrence $\chxx_{a+1} = \Tr{\xb} \chxx_a - \chxx_{a-1}$ can be derived. These polynomials, shown in Table \ref{t:rank1centralfunctions}, are closely related to the Chebyshev polynomials of the second kind.
    \begin{table}
    \begin{align*}
        \chi_0 &= 1 \\
        \chi_1 &= x \\
        \chi_2 &= x^2-1 \\
        \chi_3 &= x^3-2x \\
        \chi_4 &= x^4-3x^2+1 \\
        \chi_5 &= x^5-4x^3+3x.
    \end{align*}
    \caption{Rank 1 Central Functions.}\label{t:rank1centralfunctions}
    \end{table}
Note that the ring structure is not the usual polynomial structure of $\C[\Tr{\xb}]$.

\subsection{Example $r=2$}
The diagram is:
    $$
    \ch{c}{a}{b} = 
    \tikz[trivalent,every node/.style={basiclabel}]{
        \draw(0,.5)circle(.4)(0,.1)arc(-145:145:.7);
        \node[small matrix]at(-.4,.5){$\mathbf{X}_1$};
        \node[small matrix]at(.4,.5){$\mathbf{X}_2$};
        \node at(-.4,1){$a$};\node at(.5,.95){$b$};\node at(1.2,1.1){$c$};
    }
    $$

Recall the decomposition $$\C[\SL\times \SL]^{\SL}\approx\sum_{\substack{a,b\in\N\\
\Adm abc}}\cspan{\ch{c}{a}{b}},$$ where $\ch{c}{a}{b}$ corresponds to the image of
  $$\sum_{k=0}^c{\sf c}_k({\sf c}_k)^T\mapsto\sum_{k=0}^c\tbinom{c}{k}\bs{\sf c}{k}{k}$$
under the injection $V^*_c \otimes V_c\hookrightarrow V^*_a\otimes V^*_b \otimes V_a\otimes V_b$.

This inclusion is determined by the Clebsch-Gordan injection $\iota:V_c\hookrightarrow V_a\otimes
V_b.$ Hence, an explicit formula for $\iota$ provides a means to compute $\ch cab$ directly.

Since the general injections are determined by the rank 2 injections,  we now review their construction.

A few simple examples will motivate the construction of $\iota$.

For $k=1,2$, let $\xb_k=(x_{ij}^k)$ be $2\times2$ generic matrices, and let
\begin{align*}
  x&=\Tr{\xb_1}=x^1_{11}+x^1_{22},\\
  y&=\Tr{\xb_2}=x^2_{11}+x^2_{22},\\
  z&=\Tr{\xb_1 \xb_2^{-1}}=(x^1_{11}x^2_{22}+x^1_{22}x^2_{11})-(x^1_{12}x^2_{21}+x^1_{21}x^2_{12}).
\end{align*}

The map $\cup:V_0\hookrightarrow V_1\otimes V_1$ given by
$${\sf c}_0\mapsto{\sf a}_0\otimes {\sf b}_1-{\sf a}_1\otimes {\sf b}_0$$ is invariant.

More generally, the injection $V_0 \hookrightarrow V_a\otimes V_a$ is given by
\begin{equation*}
 \cup^a:{\sf c}_0\longmapsto\sum_{m=0}^{a}(-1)^m\tbinom{a}{m}{\sf a}_{a-m}\otimes{\sf b}_m.
\end{equation*}

Hence, $\ch000=1$ and $\ch011$ may be computed by:
\begin{align*}
  \ch011 & \mapsto \bs{\sf c}{0}{0}\\
    & \mapsto ({\sf a}^*_0\otimes {\sf b}_1^*-{\sf a}_1^*\otimes {\sf b}^*_0)\otimes({\sf a}_0\otimes {\sf b}_1-{\sf a}_1\otimes {\sf b}_0)\\
    & \mapsto (\bs{\sf a}{0}{0})\otimes(\bs{\sf b}{1}{1})-(\bs{\sf a}{0}{1})\otimes (\bs{\sf b}{1}{0}) \\
        & \hspace{.5in} -(\bs{\sf a}{1}{0})\otimes (\bs{\sf b}{0}{1})+(\bs{\sf a}{1}{1})\otimes (\bs{\sf b}{0}{0})\\
    & \mapsto x^1_{11}\otimes x^2_{22}-x^1_{12}\otimes x^2_{21}-x^1_{21}\otimes x^2_{12}+x^1_{22}\otimes x^2_{11}\\
    & \mapsto (x^1_{11}x^2_{22}+x^1_{22}x^2_{11})-(x^1_{12}x^2_{21}+x^1_{21}x^2_{12})=z.
\end{align*}

The representation $V_c$ may be identified with a subset of $\vprod c$ via the equivariant maps
  $$\xymatrix{V_c\ar@/^1pc/[r]^-{\sss\sf Sym} & \vprod{c}\ar@/^1pc/[l]^-{\sss\sf Proj}}$$
where ${\sf Proj}\circ{\sf Sym}={\rm id}$.

Thus, when $c=a+b$, $\iota$ is given by the commutative diagram
  $$\xymatrix{
     \vprod{c}\ar@{=}[r]\ar@{}[dr]|-{\text{\Large $\circlearrowright$}} & \vprod a\otimes\vprod b\ar[d]^{\sss\sf Proj\otimes Proj}\\
     V_c\ar[r]_-\iota\ar[u]^{\sss\sf Sym}              & V_a\otimes V_b.}$$

In particular,
\begin{equation*}
  \tbinom{c}{k}{\sf c}_k\overset{\iota}{\longmapsto}\sum_{\substack{0\le i \le a\\ 0 \leq j \leq b \\
  i+j=k}}\tbinom{a}{i}{\sf a}_i\otimes \tbinom{b}{j}{\sf b}_j.
\end{equation*}

For example, consider $\ch110$. In this case, ${\sf c}_0 \mapsto {\sf a}_0\otimes {\sf b}_0$ and
${\sf c}_1 \mapsto {\sf a}_1\otimes {\sf b}_0$.

Hence, {\small
\begin{align*}
  \ch110&\mapsto\bs{\sf c}00+\bs{\sf c}11\mapsto(\bs{\sf a}{0}{0})\otimes(\bs{\sf b}00)+(\bs{\sf a}11)\otimes(\bs{\sf b}00)\\
   &\mapsto x^1_{11}\otimes 1+x^1_{22}\otimes 1\mapsto x^1_{11}+x^1_{22}=x.
\end{align*}}

A similar computation shows that $\ch101\mapsto y$.

Let $\gamma=(a+b-c)/2$, $\alpha=(b+c-a)/2$, and  $\beta=(a-b+c)/2$.  The general form of $\iota$ is determined by combining these cases in the following diagram (see \eqref{eq:trivalent-vertex} for the corresponding spin network diagram):
  $$\xymatrix{
      V_c\ar[r]^-\iota\ar[d]_-\iota\ar@{}[dr]|-{\text{\Large $\circlearrowright$}}
                                      & V_{\beta}\otimes V_{\alpha}\ar[d]^{{\rm id}\otimes{\cup^\gamma}\otimes{\rm id}}\\
      V_a\otimes V_b                  & V_{\beta}\otimes V_\gamma\otimes V_\gamma\otimes V_{\alpha}\ar[l].}$$

It follows that the mapping $\iota:V_c\to V_a\otimes V_b$ is explicitly given by:
\begin{align*}
  \tbinom{c}{k}{\sf c}_k
    &\longmapsto\sum_{\substack{0\le i\le\beta\\0\le j\le\alpha\\0\le m\le\gamma\\i+j=k}}%
      \tbinom{\beta}{i}{\sf a}_i
        \otimes\left[(-1)^m\tbinom{\gamma}{m}{\sf a}_{\gamma-m}\otimes{\sf b}_m\right]
        \otimes\tbinom{\alpha}{j}{\sf b}_j\\ %
    &\longmapsto\sum_{\substack{0\le i\le\beta\\0\le j\le\alpha\\0\le m\le\gamma\\i+j=k}}%
      (-1)^m\tbinom{\beta}{i}\tbinom{\alpha}{j}\tbinom{\gamma}{m}{\sf a}_{i+\gamma-m}\otimes{\sf b}_{j+m}.%
\end{align*}

In \cite{LP} it is shown
\begin{theorem}\label{t:ranktworecurrencex}
Provided $a>1$ and $c>1$, we can write
\begin{multline*}
  \ch cab=x\cdot\ch {a-1}{b}{c-1}-\tfrac{(a+b-c)^2}{4a(a-1)}\ch c{a-2}b%
    -\tfrac{(-a+b+c)^2}{4c(c-1)}\ch {c-2}ab\\
    -\tfrac{(a+b+c)^2(a-b+c-2)^2}{16a(a-1)c(c-1)}\ch{c-2}{a-2}b.
\end{multline*}
The relation still holds for $a=1$ or $c=1$, provided we exclude the terms with $a-1$ or $c-1$ in
the denominator.
\end{theorem}

Also, note that formulae for multiplication by $y$ and $z$ may be obtained by applying the following symmetry relation.

Suppose a central function is expressed as a polynomial $P$ in the variables $x=\Tr{\xb_1}$, $y=\Tr{\xb_2}$, and
$z=\Tr{\xb_1\xb_2^{-1}}$, so that $P_{\sss a,b,c}(y,x,z)=\ch cab(\xb_1,\xb_2)$ for some admissible triple
$(a,b,c)$.

\begin{theorem}For any permutation $\sigma$,
$$P_{\sss \sigma(a,b,c)}(y,x,z)=P_{\sss a,b,c}(\sigma^{-1}(y,x,z)).$$
\end{theorem}

Using this symmetry and the above recursion, the ring structure is completely determined. As stated in the introduction, for $r>2$, this ring structure is not known and for $r=1,2$ it was worked out in \cite{LP}.%
	\footnote{The referee pointed out a typo in the statement of the relevant theorem for rank $r=2$ in \cite{LP}: the term $\{b,j_i,l\}$ in Theorem 5.11 and the corresponding terms in Lemma 5.10 should be removed from the list of admissible triples. With this correction, the theorem and its proof as  written remain true.}

Section \ref{s:rankthree} explores the $r=3$ case using computations made with {\it Mathematica}.  We use both the tensorial contraction method discussed above (which reflects our definition of central functions), and a purely combinatorial method that uses spin network techniques. The next few sections lay the groundwork for the combinatorial method, which comes from the representation of the central functions as spin networks.

%% file: tracediagramreview.tex
\subsection{Gluing Lemmas}\label{s:gluinglemmas}

Spin networks satisfy certain ``recoupling'' identities.

\begin{definition}[$6j$-symbols]\label{d:6jsymbol}
    The \emph{$6j$-symbols} are the coefficients in the following change-of-basis equation:
    $$
    \tikz[trivalent]{
        \coordinate(vxa)at(0,.35)
            edge[]node[leftlabel,pos=1]{$d$}(0,-.15)edge[bend left]node[leftlabel,pos=1]{$a$}(-.6,1.15);
        \coordinate(vxb)at(.3,.75)
            edge[bend left]node[basiclabel,pos=.7,right=2pt]{$e$}(vxa)
            edge[bend left]node[leftlabel,pos=1]{$b$}(0,1.15)edge[bend right]node[rightlabel,pos=1]{$c$}(.6,1.15);}
    = \sum_{f\in\lceil a,b\rfloor\cap\lceil c,d\rfloor} \sixj abdcef
    \tikz[trivalent]{
        \coordinate(vxa)at(0,.35)
            edge[]node[leftlabel,pos=1]{$d$}(0,-.15)edge[bend right]node[rightlabel,pos=1]{$c$}(.6,1.15);
        \coordinate(vxb)at(-.3,.75)
            edge[bend right]node[basiclabel,pos=.7,left=3pt]{$f$}(vxa)
            edge[bend right]node[rightlabel,pos=1]{$b$}(0,1.15)edge[bend left]node[leftlabel,pos=1]{$a$}(-.6,1.15);}
    .$$
\end{definition}
Note that this definition differs slightly from \cite{CFS95,Kau91,LP}.

Together with the fusion identity, this provides an identity for gluing a strand ``across'' a vertex.
\begin{proposition}[Strand-Vertex Gluing]\label{l:vertexglue}
    $$
    \tikz[trivalent]{
        \coordinate(vxa)at(.5,.5)
        edge[]node[rightlabel,pos=1]{$a$}(.5,-.4)
        edge[bend left]node[rightlabel,pos=1]{$c$}(0,1.1)
        edge[bend right]node[rightlabel,pos=1]{$b$}(1,1.1);
        \draw[blue](.3,-.4)to(.3,-.1)to[wavyup](-.2,1.1)node[leftlabel]{$d$};}
    = \sum_{e\in\iadm{a,d}} \mff_e(a,d)
    \tikz[trivalent]{
        \coordinate(vxa)at(0,.35)
            edge[blue,bend left]node[leftlabel,pos=1]{$d$}(-.6,1.15);
        \coordinate(vxc)at(0,-.15)
            edge[]node[leftlabel]{$e$}(vxa)
            edge[blue,bend right]node[leftlabel,pos=1]{$d$}(-.3,-.55)edge[bend left]node[rightlabel,pos=1]{$a$}(.3,-.55);
        \coordinate(vxb)at(.3,.75)
            edge[bend left]node[basiclabel,pos=.7,right=2pt]{$a$}(vxa)
            edge[bend left]node[leftlabel,pos=1]{$c$}(0,1.15)edge[bend right]node[rightlabel,pos=1]{$b$}(.6,1.15);}
    = \sum_{\substack{e\in\iadm{a,d}\\ f\in\iadm{c,d}\cap\iadm{b,e}}}
    \mff_e(a,d)\sixj dcebaf
    \tikz[trivalent]{
        \coordinate(vxa)at(0,.35)
            edge[bend right]node[rightlabel,pos=1]{$b$}(.6,1.15);
        \coordinate(vxc)at(0,-.15)
            edge[]node[leftlabel]{$e$}(vxa)
            edge[blue,bend right]node[leftlabel,pos=1]{$d$}(-.3,-.55)edge[bend left]node[rightlabel,pos=1]{$a$}(.3,-.55);
        \coordinate(vxb)at(-.3,.75)
            edge[bend right]node[basiclabel,pos=.7,left=3pt]{$f$}(vxa)
            edge[bend right]node[rightlabel,pos=1]{$c$}(0,1.15)edge[blue,bend left]node[leftlabel,pos=1]{$d$}(-.6,1.15);}
    .$$
\begin{proof}
    The first step applies the fusion identity (Proposition \ref{p:bubblefusion}) to the lower part of the diagram, introducing a summation over $e$. The second step uses the change-of-basis formula (Definition \ref{d:6jsymbol}) on the upper half of the diagram.
\end{proof}
\end{proposition}

We now introduce specialized notation for the case of the above lemma with $d=1$, since it will simplify the expression of the product in Theorem \ref{t:simplerecurrence}.

\begin{definition}\label{d:fusioncoeff}
    Given an admissible triple $\{a,b,c\}$, $a'\in\iadm{1,a}$, and $c'\in\iadm{1,c}$, define the \emph{fusion coefficient} $\mathfrak{F}$ and the \emph{normalized fusion coefficient} $\hat{\mathfrak{F}}$ by
    \begin{align*}
        \Fus ba{a'}c{c'} \equiv \mfs_{a'}(1,a) \mff_{a'}(1,a) \sixjt 1c{a'}ba{c'}.\\
        \nFus ba{a'}c{c'} \equiv \mfs_{a'}(1,a) \sqrt{\tfrac{\mff_{a'}(1,a)}{\mff_{c'}(1,c)}} \sixjt 1c{a'}ba{c'}.
    \end{align*}
\end{definition}

It is immediate that
    \begin{equation}\label{eq:fusnormalization}
        \Fus ba{a'}c{c'}=\sqrt{\mff_{a'}(1,a)\mff_{c'}(1,c)} \nFus ba{a'}c{c'}.
    \end{equation}
Also, $a'\in\iadm{1,a}$ is equivalent to requiring $a'=a\pm1$, so given $\{a,b,c\}$ there are four choices for $\{a',c'\}$. Table \ref{t:fusioncoeff} shows the values of the $6j$-symbol, the fusion coefficient, and the reduced coefficient in each of the four cases. One can show that for any $(a', c')\in \lceil 1, a\rfloor \times \lceil 1, c\rfloor$ that 
    $$\Fus ba{a'}c{c'}=\Fus bc{c'}a{a'} \quad\text{and}\quad \nFus ba{a'}c{c'}=\nFus bc{c'}a{a'}.$$
These values are taken from Corollary 5.5 in \cite{LP}, together with the fact that the additional fusion constant is either $\mff_{a+1}(1,a)=1$ or $\mff_{a-1}(1,a)=\frac{a}{a+1}$. This table corrects sign errors in the $6j$-symbol formulas provided on pp. 79-80 of \cite{Pet06}.

    \def\spacer{\tikz{\draw[white,transparent](0,-.2)--(0,1.2);}}
    \def\spacerb{\tikz{\draw[white,transparent](0,-.4)--(0,1.4);}}
    \begin{center}
    \begin{table}
    \begin{tabular}{|c|ccc|}
    \hline
    \spacer
    $\{a',c'\}$ & $\sixjt 1c{a'}ba{c'}$ & $\Fus ba{a'}c{c'}$ & $\nFus ba{a'}c{c'}$ \\\hline
    \spacer
    $\{a+1,c+1\}$ & $1$ & $1$ & $1$\\
    \spacer
    $\{a-1,c+1\}$ & $\tfrac{\mfe_c(a,b)}{a}$ & $-\tfrac{\mfe_c(a,b)}{a+1}$
        & $-\tfrac{\mfe_c(a,b)}{\sqrt{a(a+1)}}$ \\
    \spacer
    $\{a+1,c-1\}$ & $-\tfrac{\mfe_a(c,b)}{c+1}$ & $-\tfrac{\mfe_a(c,b)}{c+1}$
        & $-\tfrac{\mfe_a(c,b)}{\sqrt{c(c+1)}}$ \\
    \spacerb
    $\{a-1,c-1\}$ & $\frac{\mfe_b(a,c)(\mfe(a,b,c)+1)}{a(c+1)}$ & $-\tfrac{\mfe_b(a,c)(\mfe(a,b,c)+1)}{(a+1)(c+1)}$
        & $-\tfrac{\mfe_b(a,c)(\mfe(a,b,c)+1)}{\sqrt{a(a+1)c(c+1)}}$ \\\hline
    \end{tabular}
    \caption{The four basic fusion coefficients, with corresponding $6j$-symbols and normalized values.}\label{t:fusioncoeff}
    \end{table}
    \end{center}

With a different orientation, Proposition \ref{l:vertexglue} becomes
\begin{proposition}\label{l:gluevertex-alt}
    \begin{equation}\label{eq:glueonesum2}
    \tikz[trivalent]{
        \coordinate(vxa)at(.5,.4)edge[]node[rightlabel,pos=1]{$b$}(.5,-.1)
            edge[bend left]node[leftlabel,pos=1]{$c$}(0,1.1)
            edge[bend right]node[rightlabel,pos=1]{$a$}(1,1.1);
        \draw[blue](.2,1.1)to[out=-90,in=-90,looseness=2](.8,1.1)node[leftlabel]{$1$};}
    = \sum_{\substack{a'\in\iadm{1,a}\\ c'\in\iadm{1,c}\cap\iadm{b,a'}}}
        \Fus ba{a'}c{c'}
    \tikz[trivalent]{
        \coordinate(vxa)at(.5,.4)edge[]node[rightlabel,pos=1]{$b$}(.5,-.1);
        \coordinate(vxb)at(.1,.8)edge[bend right]node[leftlabel]{$c'\:$}(vxa)
            edge[bend left]node[leftlabel,pos=1]{$c$}(-.2,1.2)
            edge[blue,bend right]node[leftlabel,pos=1]{$1$}(.4,1.2);
        \coordinate(vxc)at(.9,.8)edge[bend left]node[rightlabel]{$a'$}(vxa)
            edge[blue,bend left]node[rightlabel,pos=1]{$1$}(.6,1.2)
            edge[bend right]node[rightlabel,pos=1]{$a$}(1.2,1.2);}.
    \end{equation}
\begin{proof}
    In terms of the fusion coefficient, Proposition \ref{l:vertexglue} is
    \begin{equation}\label{eq:gluevertexleft}
    \tikz[trivalent]{
        \coordinate(vxa)at(.5,.5)
        edge[]node[rightlabel,pos=1]{$a$}(.5,-.4)
        edge[bend left]node[rightlabel,pos=1]{$c$}(0,1.1)
        edge[bend right]node[rightlabel,pos=1]{$b$}(1,1.1);
        \draw[blue](.3,-.4)to(.3,-.1)to[wavyup](-.2,1.1)node[leftlabel]{$1$};}
    = \sum_{\substack{a'\in\iadm{a,1}\\ c'\in\iadm{c,1}\cap\iadm{b,a'}}}
    \mfs_{a'}(a,1)\Fus ba{a'}c{c'}
    \tikz[trivalent]{
        \coordinate(vxa)at(0,.35)
            edge[bend right]node[rightlabel,pos=1]{$b$}(.6,1.15);
        \coordinate(vxc)at(0,-.15)
            edge[]node[leftlabel]{$a'$}(vxa)
            edge[bend right]node[blue,leftlabel,pos=1]{$1$}(-.3,-.55)edge[bend left]node[rightlabel,pos=1]{$a$}(.3,-.55);
        \coordinate(vxb)at(-.3,.75)
            edge[bend right]node[basiclabel,pos=.7,left=3pt]{$c'$}(vxa)
            edge[bend right]node[rightlabel,pos=1]{$c$}(0,1.15)edge[blue,bend left]node[leftlabel,pos=1]{$1$}(-.6,1.15);}.
    \end{equation}
    Reflect this relation vertically using Proposition \ref{p:spinnetreflection}, and extend the strands labeled by $1$ and $c$ on both sides of the equation to obtain:
        \begin{equation}\label{eq:extendstrands}
        \tikz[trivalent]{
            \coordinate(vertex)at(.5,.5)edge[]node[rightlabel,pos=1]{$a$}(.5,1)
                edge[bend right](.1,.2)edge[bend left]node[rightlabel,pos=1]{$b$}(1,0);
            \draw[blue](-.1,.2)to[wavyup](.3,1);
            \draw[red](-.1,.2)to[out=-90,in=-90](-.5,.2)to[wavyup](-.4,1)node[rightlabel]{$1$}
                (.1,.2)to[out=-90,in=-90,looseness=1.5](-.7,.2)to[wavyup](-.6,1)node[leftlabel]{$c$};
        }=\sum_{\substack{a'\in\iadm{a,1}\\ c'\in\iadm{c,1}\cap\iadm{b,a'}}}
            \mfs_{a'}(a,1)\Fus ba{a'}c{c'}
        \tikz[trivalent]{
            \coordinate(vertexa)at(.2,.2)edge[blue,bend right](-.1,0)edge[bend left](.3,0);
            \coordinate(vertexb)at(.5,.5)edge[bend right]node[leftlabel,pos=.3]{$c'$}(vertexa)
                edge[bend left]node[rightlabel,pos=1]{$b$}(1.1,-.2);
            \coordinate(vertexc)at(.5,.8)edge[]node[rightlabel,pos=.3]{$a'$}(vertexb)
                edge[blue,bend left]node[leftlabel,pos=1]{$1$}(.3,1.2)
                edge[bend right]node[rightlabel,pos=1]{$a$}(.7,1.2);
            \draw[red](-.1,0)to[out=-90,in=-90](-.5,0)to[wavyup](-.6,1.2)node[rightlabel]{$1$}
                (.3,0)to[out=-90,in=-90,looseness=1.5](-.7,0)to[wavyup](-.8,1.2)node[leftlabel]{$c$};
        }.
        \end{equation}
    By Proposition \ref{p:spinnetsignstrong}, the relation can be straightened, with the introduction of signs $S_1$ and $S_2$, to:
        \begin{equation}\label{eq:straightenstrands}
        S_1
        \tikz[trivalent]{
            \coordinate(vxa)at(.5,.4)edge[]node[rightlabel,pos=1]{$b$}(.5,-.1)
                edge[bend left]node[leftlabel,pos=1]{$c$}(0,1.1)
                edge[bend right]node[rightlabel,pos=1]{$a$}(1,1.1);
            \draw[blue](.2,1.1)to[out=-90,in=-90,looseness=2](.8,1.1)node[leftlabel]{$1$};}
        = \sum_{\substack{a'\in\iadm{a,1}\\ c'\in\iadm{c,1}\cap\iadm{b,a'}}}
            S_2\:
            \mfs_{a'}(a,1)\Fus ba{a'}c{c'}
        \tikz[trivalent]{
            \coordinate(vxa)at(.5,.4)edge[]node[rightlabel,pos=1]{$b$}(.5,-.1);
            \coordinate(vxb)at(.1,.8)edge[bend right]node[leftlabel,pos=.5]{$c'\:$}(vxa)
                edge[bend left]node[leftlabel,pos=1]{$c$}(-.2,1.2)
                edge[blue,bend right]node[leftlabel,pos=1]{$1$}(.4,1.2);
            \coordinate(vxc)at(.9,.8)edge[bend left]node[rightlabel,pos=.5]{$a'$}(vxa)
                edge[blue,bend left]node[rightlabel,pos=1]{$1$}(.6,1.2)
                edge[bend right]node[rightlabel,pos=1]{$a$}(1.2,1.2);}.
        \end{equation}

To calculate $S_1S_2$, we must count the number of kinks in each diagram. No strands in \eqref{eq:straightenstrands} are kinked. The strands between $b$ and $c$ on the left-hand side of \eqref{eq:extendstrands} are kinked, giving a factor $S_1=\mfs_a(b,c)$. On the right-hand side, strands between $b$ and either $1$ or $c$ are kinked, and so are those between $1$ and $c$, producing a factor $S_2=\mfs_{a'}(b,c')\mfs_{c'}(1,c)$. Hence
        $$S_1S_2=\mfs_{c'}(1,c)\mfs_{a'}(b,c')\mfs_a(b,c)=(-1)^{\frac12(1-(a'-a))}=\mfs_{a'}(1,a).$$
    Therefore, the coefficient is $\left(\mfs_{a'}(1,a)\right)^2\Fus ba{a'}c{c'}=\Fus ba{a'}c{c'}$.
\end{proof}
\end{proposition}

The reflective symmetries guaranteed by Proposition \ref{p:spinnetreflection} imply that corresponding formulas for the products
    \begin{equation}\label{eq:gluevertexlist}
    \tikz[trivalent]{
        \coordinate(vxa)at(.5,.5)
            edge[]node[leftlabel,pos=1]{$a$}(.5,-.1)
            edge[bend left]node[leftlabel,pos=1]{$b$}(0,1.1)
            edge[bend right]node[leftlabel,pos=1]{$c$}(1,1.1);
        \draw[blue](.7,-.1)to[wavyup](1.2,1.1)node[rightlabel]{$1$};}
    \quad
    \tikz[trivalent]{
        \coordinate(vxa)at(.5,.5)
            edge[]node[rightlabel,pos=1]{$a$}(.5,1.1)
            edge[bend left]node[rightlabel,pos=1]{$b$}(1,-.1)
            edge[bend right]node[rightlabel,pos=1]{$c$}(0,-.1);
        \draw[blue](-.2,-.1)to[wavyup](.3,1.1)node[leftlabel]{$1$};}
    \quad
    \tikz[trivalent]{
        \coordinate(vxa)at(.5,.5)
            edge[]node[leftlabel,pos=1]{$a$}(.5,1.1)
            edge[bend left]node[leftlabel,pos=1]{$c$}(1,-.1)
            edge[bend right]node[leftlabel,pos=1]{$b$}(0,-.1);
        \draw[blue](1.2,-.1)to[wavyup](.7,1.1)node[rightlabel]{$1$};}
    \end{equation}
can be generated by reflecting \eqref{eq:gluevertexleft}. In particular, the terms in these summations have the same coefficients.


%% file: simplerecurrence.tex
\subsection{Simple Loop Recurrences}\label{simplerecurrence-section}\label{s:simple}

Relation \eqref{eq:gluevertexleft} can be ``stacked'' to obtain a more general formula for gluing across two or more vertices.

\begin{lemma}\label{l:vertexglue2}
    If the edges labeled by $\{a_0,\ldots,a_n\}$ are disjoint, then
        $$
        \tikz[trivalent,shift={(0,-2)}]{
            \foreach\xa/\xb/\xc in{0/0/a_0,-.1/.8/a_1,-.4/3/a_{n-1},-.5/3.8/a_n}{
                \draw(\xa,\xb)to[out=110,in=-90]([shift={(\xa,\xb)}]-.1,.8)node[rightlabel,pos=.8]{$\xc$};
                \draw[blue,shift={(-.35,0)}](\xa,\xb)to[wavyup]([shift={(\xa,\xb)}]-.1,.8);
            }
            \foreach\xa/\xb/\xc in{-.1/.8/b_1,-.5/3.8/b_n}{
                \draw(\xa,\xb)to[bend right]([shift={(\xa,\xb)}].6,.7)node[rightlabel]{$\xc$};
            }
            \draw[dotted](-.2,1.6)to[wavyup](-.4,3);
            \draw[blue,dotted,shift={(-.35,0)}](-.2,1.6)to[wavyup](-.4,3);
            \node[blue,leftlabel]at(-1,4.6){$1$};
        }
        = \sum_{\substack{a_0'\in\iadm{1,a_0}\\ \{a_i'\in\iadm{1,a_i}\cap\iadm{b_i,a_{i-1}'}\}_{i=1}^n}}
   \hspace{-1.4cm}\sqrt{\mff_{a_0'}(a_0,1)\mff_{a_n'}(a_n,1)}
            \prod_{i=1}^n\mfs_{a_{i-1}'}(a_{i-1},1)\nFus{b_i}{a_{i-1}}{a_{i-1}'}{a_i}{a_i'}
        \tikz[trivalent,shift={(0,-2)}]{
            \foreach\xa/\xb/\xc in{-.1/.8/a_1',-.4/3/a_{n-1}'}{
                \draw(\xa,\xb)to[out=110,in=-90]([shift={(\xa,\xb)}]-.1,.8)node[rightlabel,pos=.8]{$\xc$};}
            \foreach\xa/\xb/\xc in{-.1/.8/b_1,-.5/3.8/b_n}{
                \draw(\xa,\xb)to[bend right]([shift={(\xa,\xb)}].6,.7)node[rightlabel]{$\xc$};}
            \draw[dotted](-.2,1.6)to[wavyup](-.4,3);
            \coordinate(vxa)at(0,.3)
                edge[wavyup]node[rightlabel]{$a_0'$}(-.1,.8)
                edge[blue,bend right]node[leftlabel,pos=1]{$1$}(-.3,-.1)
                edge[bend left]node[rightlabel,pos=1]{$a_0$}(.2,-.1);
            \coordinate(vxb)at(-.6,4.3)
                edge[wavydown]node[leftlabel]{$a_n'$}(-.5,3.8)
                edge[blue,bend left]node[leftlabel,pos=1]{$1$}(-.9,4.7)
                edge[bend right]node[rightlabel,pos=1]{$a_n$}(-.4,4.7);
        }.
        $$
\begin{proof}
    For clarity, we present the concrete $n=2$ case here, from which the general pattern can be seen. First, apply \eqref{eq:gluevertexleft} on the upper and lower halves of the diagram separately:
    \begin{multline*}
    \tikz[trivalent,shift={(0,-.3)}]{
        \draw[very thick,lightgray](-.15,.5)to(1.45,.5);
        \foreach\xx/\xy/\xa/\xb/\xcc/\xd in{.25/-.6/a_0/b_1//,-.25/.6/a_1/b_2/a_2/1}{\begin{scope}[shift={(\xx,\xy)}]
        \coordinate(vxa)at(.5,.5)
            edge[]node[rightlabel,pos=.45]{$\xa$}(.5,-.1)
            edge[bend left]node[rightlabel,pos=1]{$\xcc$}(0,1.1)
            edge[bend right]node[rightlabel,pos=.7]{$\xb$}(1,1.1);
        \draw[blue](.3,-.1)to[wavyup](-.2,1.1)node[leftlabel]{$\xd$};
        \end{scope}}
    }
    = \sum_{\substack{a_0'\in\iadm{1,a_0}\\
        a_1'\in\iadm{1,a_1}\cap\iadm{b_1,a_0'}\\
        a_2'\in\iadm{1,a_2}\cap\iadm{b_2,a_1'}}}
    \mfs_{a_0'}(a_0,1)\Fus{b_1}{a_0}{a_0'}{a_1}{a_1'}
    \mfs_{a_1'}(a_1,1)\Fus{b_2}{a_1}{a_1'}{a_2}{a_2'}
    \tikz[trivalent]{
        \draw[very thick,lightgray](-.65,.3)to(.95,.3);
        \foreach\xx/\xy/\xa/\xb/\xcc/\xd in{.15/-.85/a_0/b_1/a_1/1,-.15/.85/a_1/b_2/a_2/1}
            {\begin{scope}[shift={(\xx,\xy)}]
        \coordinate(vxa)at(0,.35)
            edge[bend right]node[rightlabel,pos=.7]{$\xb$}(.6,1.15);
        \coordinate(vxc)at(0,-.15)
            edge[]node[leftlabel]{$\xa'$}(vxa)
            edge[blue,bend right]node[leftlabel,pos=1]{$\xd$}(-.3,-.55)
            edge[bend left]node[rightlabel,pos=1]{$\xa$}(.3,-.55);
        \coordinate(vxb)at(-.3,.75)
            edge[bend right]node[pos=.7,pin={[basiclabel]left:$\xcc'$}]{}(vxa)
            edge[bend right]node[rightlabel,pos=1]{$\xcc$}(0,1.15)
            edge[blue,bend left]node[leftlabel,pos=1]{$\xd$}(-.6,1.15);
        \end{scope}}}
    \\
    = \sum_{\substack{a_0'\in\iadm{1,a_0}\\ \{a_i'\in\iadm{1,a_i}\cap\iadm{b_i,a_{i-1}'}\}_{i=1}^2}}
    \mfs_{a_0'}(a_0,1)\Fus{b_1}{a_0}{a_0'}{a_1}{a_1'}
    \mfs_{a_1'}(a_1,1)\Fus{b_2}{a_1}{a_1'}{a_2}{a_2'} \mfb_{a_1'}(1,a_1)
    \tikz[trivalent,scale=.9,shift={(0,-.2)}]{
        \coordinate(vxa)at(0,.35)
            edge[bend right]node[rightlabel,pos=.7]{$b_1$}(.6,1.15)
            edge[bend left](-.3,1.15);
        \coordinate(vxc)at(0,-.15)
            edge[]node[leftlabel]{$a_0'$}(vxa)
            edge[blue,bend right]node[leftlabel,pos=1]{$1$}(-.3,-.55)
            edge[bend left]node[rightlabel,pos=1]{$a_0$}(.3,-.55);
        \begin{scope}[shift={(-.3,1.7)}]
        \coordinate(vxa)at(0,.35)
            edge[bend right]node[rightlabel,pos=.7]{$b_2$}(.6,1.15)
            edge[]node[leftlabel]{$a_1'$}(0,-.7);
        \coordinate(vxb)at(-.3,.75)
            edge[bend right]node[pos=.7,pin={[basiclabel]left:$a_2'$}]{}(vxa)
            edge[bend right]node[rightlabel,pos=1]{$a_2$}(0,1.15)
            edge[blue,bend left]node[leftlabel,pos=1]{$1$}(-.6,1.15);
        \end{scope}
    }
    \end{multline*}
Note that the bubble identity \eqref{eq:bubbleidentity} allows us to use the same index $a_1'$ above and below the bubble in the second diagram and contributes the $\mfb_{a_1'}(1,a_1)$ factor in the second step. Using \eqref{eq:fusnormalization} and the fact that $\mff_c(a,b)\mfb_c(a,b)=1$, the coefficient can be expressed as
    $$\sqrt{\mff_{a_0'}(a_0,1)\mff_{a_2'}(a_2,1)}
        \prod_{i=1}^2\mfs_{a_{i-1}'}(a_{i-1},1)\nFus{b_i}{a_{i-1}}{a_{i-1}'}{a_i}{a_i'}.
    $$
    The case $n>2$ follows by induction.
\end{proof}
\end{lemma}


\begin{definition}
		Let $\mathsf{s}$ be a trivalent graph with edges $E=\{e_1,e_2,\ldots,e_m\}$. An \emph{admissible labeling} of the edges of $\mathsf{s}$ is a map $l:E\to\mathbb{N}$ such that the labels at all vertices form admissible triples. The resulting spin network is denoted $\mathsf{s}_l$.

    An \emph{admissible relabeling} of $l$ is a labeling of a subset of edges $a:\{e_{j_1},e_{j_2},\ldots,e_{j_n}\}\to\mathbb{N}$ such that the labeling defined by
    $$l_a(e_i)=\left\{
    	\begin{array}{ll}
    		a(e_i) & \text{if } i=j_k \text{ for some $k$} \\
    		l(e_i) & otherwise.
    	\end{array}
    \right.$$
is also admissible.
\end{definition}

\begin{definition}
A \emph{simple cycle} $(e_{1},e_{2},\ldots,e_{n})$ is an ordered tuple of distinct edges in a graph that begins and ends at the same vertex.
\end{definition}

The next theorem says that multiplication of a spin network by a loop drawn parallel to a simple cycle may be expanded in terms of diagrams with admissible relabelings.

\begin{theorem}[Simple Loop Multiplication Formula]\label{t:simplerecurrence}
    Let $\mathsf{s}$ be a trivalent graph with admissible labeling $$l:\{e_1,e_2,\ldots,e_m\}\to\mathbb{N}.$$ Let $(e_{j_1},e_{j_2},\ldots,e_{j_n})$ be a simple cycle in $\mathsf{s}$, and let $\gamma$ denote the spin network consisting of a single loop drawn parallel to this cycle and labeled by 1.
    
    Define $a_i=l(e_{j_i})$, and let $b_i\in\iadm{a_{i-1},a_i}$ be the third label on the vertex joining edges $e_{j_{i-1}}$ and $e_{j_i}$. As illustrated in Figure \ref{f:ni-definition}, let $N_i$ be the sum of (a) the number of times $\gamma$ crosses the edge $e_{j_i}$, (b) the number of local extrema along the edge $e_{j_i}$ that do not occur at a vertex, and (c) the number of times $e_{j_i}$ adjoins a vertex. Then
    $$\gamma\cdot\mathsf{s}_l
     =\sum_{a'}
        \left(\prod_{i=1}^n\mfs_{a_i'}(1,a_i)^{N_i}\nFus{b_i}{a_{i-1}}{a_{i-1}'}{a_i}{a_i'}\right)
        \mathsf{s}_{l_{a'}},$$
	where the summation is over all admissible relabelings $$a':\{e_{j_1},e_{j_2},\ldots,e_{j_n}\}\to\N$$ satisfying $a_i'\in\iadm{1,a_i}\cap\iadm{b_i,a_{i-1}'}$ for all $i$, with the understanding that $e_{j_0}=e_{j_n}$ and $a_0=a_n$.
    \begin{figure}[htb]
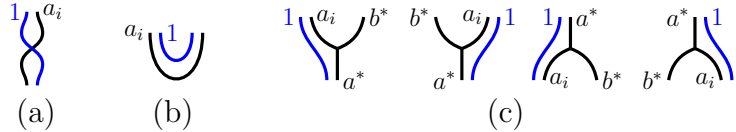

    \begin{tabular}{ccccc}
    \tikz[trivalent]{
        \draw(-.1,-.2)to[wavyup](-.2,.2)to[wavyup](.2,.8)to[wavyup](.1,1.2)node[rightlabel]{$a_i$};
        \draw[blue](.1,-.2)to[wavyup](.2,.2)to[wavyup](-.2,.8)to[wavyup](-.1,1.2)node[leftlabel]{$1$};
    }
    & &
    \tikz[trivalent]{
        \draw(-.5,.8)node[leftlabel]{$a_i$}to[out=-90,in=-90,looseness=3](.5,.8);
        \draw[blue](-.3,.8)node[rightlabel]{$1$}to[out=-90,in=-90,looseness=3](.3,.8);
    }\:\:\:
    & &
    \tikz[trivalent]{
        \coordinate(vxa)at(.5,.5)
            edge[]node[rightlabel,pos=1]{$a^*$}(.5,-.1)
            edge[bend left]node[rightlabel,pos=1]{$a_i$}(0,1.1)
            edge[bend right]node[rightlabel,pos=1]{$b^*$}(1,1.1);
        \draw[blue](.3,-.1)to[wavyup](-.2,1.1)node[leftlabel]{$1$};}
    \tikz[trivalent]{
        \coordinate(vxa)at(.5,.5)
            edge[]node[leftlabel,pos=1]{$a^*$}(.5,-.1)
            edge[bend left]node[leftlabel,pos=1]{$b^*$}(0,1.1)
            edge[bend right]node[leftlabel,pos=1]{$a_i$}(1,1.1);
        \draw[blue](.7,-.1)to[wavyup](1.2,1.1)node[rightlabel]{$1$};}
    \tikz[trivalent]{
        \coordinate(vxa)at(.5,.5)
            edge[]node[rightlabel,pos=1]{$a^*$}(.5,1.1)
            edge[bend left]node[rightlabel,pos=1]{$b^*$}(1,-.1)
            edge[bend right]node[rightlabel,pos=1]{$a_i$}(0,-.1);
        \draw[blue](-.2,-.1)to[wavyup](.3,1.1)node[leftlabel]{$1$};}
    \tikz[trivalent]{
        \coordinate(vxa)at(.5,.5)
            edge[]node[leftlabel,pos=1]{$a^*$}(.5,1.1)
            edge[bend left]node[leftlabel,pos=1]{$a_i$}(1,-.1)
            edge[bend right]node[leftlabel,pos=1]{$b^*$}(0,-.1);
        \draw[blue](1.2,-.1)to[wavyup](.7,1.1)node[rightlabel]{$1$};
    }
    \\
    {\rm(a)} & & {\rm(b)} & & {\rm(c)}
    \end{tabular}
    \caption{Local contributors to sign of simple loop multiplication formula, where either $a^*=a_{i-1}$ and $b^*=b_i$ or $a^*=a_{i+1}$ and $b^*=b_{i+1}$. The arc labeled by $1$ is part of a simple loop $\gamma$ drawn parallel to a simple cycle.}\label{f:ni-definition}
    \end{figure}
\begin{proof}
    First, suppose the loop has the form of Lemma \ref{l:vertexglue2}, with the $a_0$ and $a_n$ edges coinciding, as follows:
    $$
        \tikz[trivalent,shift={(0,.2)}]{
            \draw[blue](-.2,.8)arc(0:180:.35)--(-.9,-.8)arc(180:360:.35)--cycle;
            \draw(0,.8)arc(0:180:.55)--(-1.1,-.8)arc(180:360:.55)--cycle;
            \draw(0,-.4)to[bend right](.5,0);\draw(0,.4)to[bend right](.5,.8);
            \node at(.25,-.1){.};\node at(.25,.1){.};\node at(.25,.3){.};
        }
        \longrightarrow
        \tikz[trivalent,shift={(0,.2)}]{
            \draw[blue](0,.6)to[bend left](-.2,.8)arc(0:180:.35)--(-.9,-.8)arc(180:360:.35)to[bend left](0,-.6);
            \draw(0,.8)arc(0:180:.55)--(-1.1,-.8)arc(180:360:.55)--cycle;
            \draw(0,-.4)to[bend right](.5,0);\draw(0,.4)to[bend right](.5,.8);
            \node at(.25,-.1){.};\node at(.25,.1){.};\node at(.25,.3){.};
        }
        \longrightarrow
        \tikz[trivalent,shift={(0,.2)}]{
            \draw(0,.8)arc(0:180:.55)--(-1.1,-.8)arc(180:360:.55)--cycle;
            \draw(0,-.4)to[bend right](.5,0);\draw(0,.4)to[bend right](.5,.8);
            \node at(.25,-.1){.};\node at(.25,.1){.};\node at(.25,.3){.};
        }
    $$
Popping the final bubble introduced along these edges results in an additional factor of $\mfb_{a_0'}(a_0,1)$, which cancels with  $\sqrt{\mff_{a_0'}(a_0,1)\mff_{a_n'}(a_n,1)}=\mff_{a_0'}(a_0,1)$; the last equality holding since $a_n=a_0$. So in this case the  coefficients of the summation are
        $$\prod_{i=1}^n \mfs_{a_i'}(1,a_i) \nFus{b_i}{a_{i-1}}{a_{i-1}'}{a_i}{a_i'}.$$

In general, the relative positions of the labels $\{a_{i-1},a_i,b_i\}$ in the vicinity of a vertex matters. If the product appears locally at a local extrema, as in \eqref{eq:glueonesum2} or its reflection, then the factor $\mfs_{a_i'}(1,a_i)$ is unnecessary. Otherwise, in the cases depicted in Figure \ref{f:ni-definition}(c), the factor remains.

If a crossing or local extrema occurs along an edge, then the ``bubble popping'' step in Lemma \ref{l:vertexglue2} becomes one of the following:
    \begin{align*}
    \tikz[trivalent]{
        \draw(-.15,-.5)to[wavyup](-.1,-.2)to[wavyup](-.2,.2)to[wavyup](.2,.8)to[wavyup](.1,1.2)to[wavyup](.15,1.5)node[rightlabel]{$a_i$};
        \draw[blue](.15,-.5)to[wavyup](.1,-.2)to[wavyup](.2,.2)to[wavyup](-.2,.8)to[wavyup](-.1,1.2)to[wavyup](-.15,1.5)node[leftlabel]{$1$};}
    :\qquad
    \tikz[trivalent]{
        \coordinate(vxa)at(0,-.4)edge[bend right]node[leftlabel,pos=1]{$a_i$}(-.3,-.6)edge[blue,bend left]node[rightlabel,pos=1]{$1$}(.3,-.6);
        \coordinate(vxb)at(0,.15)edge[]node[rightlabel]{$a_i'$}(vxa)edge[bend left](-.3,.25)edge[blue,bend right](.3,.25);
        \coordinate(vxd)at(0,1.4)edge[blue,bend left]node[leftlabel,pos=1]{$1$}(-.3,1.6)edge[bend right]node[rightlabel,pos=1]{$a_i$}(.3,1.6);
        \coordinate(vxc)at(0,.85)edge[]node[rightlabel]{$b_i'$}(vxd)edge[bend left](.3,.75)edge[blue,bend right](-.3,.75);
        \draw(-.3,.25)node[leftlabel]{$a_i$}to[wavyup](.3,.75);
        \draw[blue](.3,.25)node[rightlabel]{$1$}to[wavyup](-.3,.75);
    }
    &=\mfs_{a_i'}(1,a_i)\mfb_{a_i'}(1,a_i)
    \tikz[trivalent]{
        \coordinate(vxa)at(0,.25)
            edge[bend right]node[leftlabel,pos=1]{$a_i$}(-.3,-.2)
            edge[blue,bend left]node[rightlabel,pos=1]{$1$}(.3,-.2);
        \coordinate(vxb)at(0,.75)
            edge[]node[auto,basiclabel]{$a_i'$}(vxa)
            edge[blue,bend left]node[leftlabel,pos=1]{$1$}(-.3,1.2)
            edge[bend right]node[rightlabel,pos=1]{$a_i$}(.3,1.2);};\\
    \tikz[trivalent]{
        \draw(-.5,.8)node[leftlabel]{$a_i$}to[out=-90,in=-90,looseness=3](.5,.8);
        \draw[blue](-.3,.8)node[rightlabel]{$1$}to[out=-90,in=-90,looseness=3](.3,.8);
    }
    :\qquad
    \tikz[trivalent]{
        \coordinate(vxa)at(-.5,.7)edge[bend left]node[leftlabel,pos=1]{$a_i$}(-.8,1.1)edge[blue,bend right]node[leftlabel,pos=1]{$1$}(-.2,1.1);
        \coordinate(vxc)at(-.5,.3)edge[]node[leftlabel]{$a_i'$}(vxa)edge[bend right]node[leftlabel,pos=1]{$a_i$}(-.8,0)edge[bend left,blue](-.2,0);
        \coordinate(vxb)at(.5,.7)edge[bend right]node[rightlabel,pos=1]{$a_i$}(.8,1.1)edge[blue,bend left]node[rightlabel,pos=1]{$1$}(.2,1.1);
        \coordinate(vxd)at(.5,.3)edge[]node[rightlabel]{$a_i'$}(vxb)edge[bend left](.8,0)edge[bend right,blue](.2,0);
        \draw(-.8,0)to[out=-90,in=-90](.8,0);
        \draw[blue](-.2,0)to[out=-90,in=-90]node[basiclabel,auto]{$1$}(.2,0);
    }
    &=\mfs_{a_i'}(1,a_i)\mfb_{a_i'}(1,a_i)
    \tikz[trivalent]{
        \coordinate(vxa)at(-.5,.7)edge[bend left]node[leftlabel,pos=1]{$a_i$}(-.8,1.1)edge[blue,bend right]node[leftlabel,pos=1]{$1$}(-.2,1.1);
        \coordinate(vxb)at(.5,.7)edge[bend right]node[rightlabel,pos=1]{$a_i$}(.8,1.1)edge[blue,bend left]node[rightlabel,pos=1]{$1$}(.2,1.1);
        \draw(vxa)to[out=-90,in=-90,looseness=3]node[basiclabel,auto]{$a_i'$}(vxb);
    }.
    \end{align*}
    The signs are calculated using Lemma \ref{p:spinnetsigns} and the stronger Proposition \ref{p:spinnetsignstrong}. In the second case, the sign is calculated by comparing the $\mfe_{a_i'}(1,a_i)$ kinks in the diagram
    \tikz[trivalent,scale=.75,shift={(0,.4)}]{
        \coordinate(vxc)at(-.5,.3)
            edge[]node[leftlabel]{$a_i'$}(-.5,.8)
            edge[bend right](-.7,.1)
            edge[bend left,blue](-.2,.1);
        \draw(-.7,.1)arc(180:360:.7)--(.7,.8)node[rightlabel,auto]{$a_i$};
        \draw[blue](-.2,.1)arc(180:360:.25)--(.3,.8)node[leftlabel,auto]{$1$};
    }
    with the diagram
    \tikz[trivalent,scale=.75,shift={(0,.3)}]{
        \draw(-.5,.8)node[leftlabel]{$a_i'$}--(-.5,.2)arc(180:360:.5);
        \draw(.5,.2)to[bend left,blue](.3,.8)node[leftlabel,auto]{$1$};
        \draw(.5,.2)to[bend right](.7,.8)node[rightlabel,auto]{$a_i$};
    },
    which has no kinks.
    In each case the additional sign $\mfs_{a_i'}(1,a_i)$ adds one to the exponent $N_i$.
\end{proof}
\end{theorem}

Rearranging the terms in the above theorem provides a recurrence formula, in which each diagram can be written in terms of diagrams of lower order, where order is defined as follows.
\begin{definition}
The \emph{order} of a spin network $\mathsf{s}_l$ is the sum $\sum_{i=1}^m l(e_i)$.
\end{definition}
\begin{corollary}\label{c:simplerecurrence}
    Let $\mathsf{s}$, $l$, $\gamma$, $a_i$, $b_i$, and $N_i$ be defined as in Theorem \ref{t:simplerecurrence}. If $\na:\{e_{j_1},\ldots,e_{j_m}\}\to\N$ defined by $\na(e_{j_i})=l(e_{j_i})-1$ is an admissible relabeling of $l$, then $\mathsf{s}_l$ can be expressed in terms of $\gamma$ and spin networks of lower order:
    \begin{equation}\label{eq:simplerecurrence}
        \mathsf{s}_l = \gamma\cdot\mathsf{s}_{l_{\na}}
        - \sum_{a'}
        \left(\prod_{i=1}^n\mfs_{a_i'}(1,a_i-1)^{N_i}\nFus{b_i}{a_{i-1}-1}{a_{i-1}'}{a_i-1}{a_i'}\right)
        \mathsf{s}_{l_{a'}},
    \end{equation}
where the summation is over admissible relabelings $a':\{e_{j_1},\ldots,e_{j_m}\}\to\N$ in which $a_i'\in\iadm{1,a_i-1}\cap\iadm{b_i,a_{i-1}'}$ for all $i$ and $a_i'\ne a_i$ for some $i$.
\begin{proof}
		Apply Theorem \ref{t:simplerecurrence} to the product $\gamma\cdot\mathsf{s}_{l_{\na}}$. The original network $\mathsf{s}_l$ is the summand in which $a'(e_{j_i})=a(e_{j_i})=l(e_{j_i})$ for all $i$, and is the unique one with highest order. Its coefficient is 1 since $\mfs_{a_i}(1,a_i-1)=+1$ and $\nFus{b_i}{a_{i-1}-1}{a_{i-1}}{a_i-1}{a_i}=1$. The remaining summands comprise the summation in \eqref{eq:simplerecurrence}.
\end{proof}
\end{corollary}

\subsubsection{Some Examples}

We illustrate the application of these theorems in a few basic examples. Note that the sign is only a factor when $a_i'=a_i-1$ since
    $$\mfs_{a_i'}(1,a_i)=
    \begin{cases}
        +1  &   a_i'=a_i+1  \\
        -1  &   a_i'=a_i-1.
    \end{cases}$$
We will give the recurrences in their most general form; for some choices of labels, the non-admissible terms should be excluded.

\begin{example}\label{ex:loopone}
    A single edge loop has a two-term recurrence:
    \begin{equation}\label{eq:loopone}
    \begin{matrix}
    \tikz[trivalent]{
        \draw[blue](0,.5)circle(.4);\node[blue,basiclabel]at(.1,.5){$\gamma$};
        \draw(0,.5)circle(.6);\node[basiclabel]at(.7,0){$a$};
        \draw(0,1.1)to[bend left](.1,1.5)node[rightlabel]{$b$};}
     &=& \nFus ba{a+1}a{a+1}
     &\tikz[trivalent]{\draw(0,.5)circle(.5);\node[basiclabel]at(.7,-.1){$a+1$};\draw(0,1)to[bend left](.1,1.4)node[rightlabel]{$b$};}
     &-& \nFus ba{a-1}a{a-1}
     &\tikz[trivalent]{\draw(0,.5)circle(.5);\node[basiclabel]at(.7,-.1){$a-1$};\draw(0,1)to[bend left](.1,1.4)node[rightlabel]{$b$};}
     \\
     &=&
     &\tikz[trivalent]{\draw(0,.5)circle(.5);\node[basiclabel]at(.7,-.1){$a+1$};\draw(0,1)to[bend left](.1,1.4)node[rightlabel]{$b$};}
     &+&\tfrac{(a-\frac b2)(a+\frac b2+1)}{a(a+1)}
     &\tikz[trivalent]{\draw(0,.5)circle(.5);\node[basiclabel]at(.7,-.1){$a-1$};\draw(0,1)to[bend left](.1,1.4)node[rightlabel]{$b$};}.
    \end{matrix}
    \end{equation}
    The loop $\gamma$ consists of a single edge, with $N_1=1$ since there are no crossings and $\gamma$ has one extremum that does not occur at a vertex.
\end{example}
Setting $b=0$, one obtains a formula equivalent to \eqref{eq:rank1recurrence}, which can be used to compute the rank one central functions.

\begin{example}\label{ex:looptwo}
    A two-edge loop has the following four-term recurrence:
    \begin{equation}\label{eq:looptwo}
    \begin{matrix}
    \tikz[trivalent]{
        \draw[blue](0,.5)circle(.4);\node[blue,basiclabel]at(.1,.5){$\gamma$};
        \draw(0,.5)circle(.6);\node[basiclabel]at(-.85,.8){$a_1$};\node[basiclabel]at(.85,.2){$a_2$};
        \draw(0,1.1)to[bend left](.1,1.5)node[rightlabel]{$b_2$};
        \draw(0,-.1)to[bend left](-.1,-.5)node[leftlabel]{$b_1$};}
     &=&
     \displaystyle{\sum_{\substack{a_1'=a_1\pm1\\ a_2'=a_2\pm1}}
        \left(\prod_{i=1}^2 \nFus{b_i}{a_1}{a_1'}{a_2}{a_2'}\right)}
     \tikz[trivalent]{
        \draw(0,.5)circle(.5);\node[basiclabel]at(-.8,.7){$a_1'$};\node[basiclabel]at(.8,.3){$a_2'$};
        \draw(0,1)to[bend left](.1,1.4)node[rightlabel]{$b_2$};
        \draw(0,0)to[bend left](-.1,-.4)node[leftlabel]{$b_1$};}
    \end{matrix}
    \end{equation}
    Since $N_1=0+0+0=N_2$, no additional signs are necessary.
\end{example}

\subsubsection{Application to Rank Two Central Functions}
The rank two central function is
    $$
    \ch cab (\xb_1,\xb_2) =
    \tikz[trivalent,every node/.style={basiclabel}]{
        \draw(0,.5)circle(.4)(0,.1)arc(-145:145:.7);
        \node[small matrix]at(-.4,.5){$\mathbf{X}_1$};
        \node[small matrix]at(.4,.5){$\mathbf{X}_2$};
        \node at(-.4,1){$a$};\node at(.5,.95){$b$};\node at(1.2,1.1){$c$};
    }
    $$
There are three simple loops: $(a,b)$ corresponding to $\Tr{\xb_1\xb_2^{-1}}$, $(a,c)$ corresponding to $\Tr{\xb_1}$, and $(b,c)$ corresponding to $\Tr{\xb_2}$. Each loop provides a different recurrence, a fact which was used in \cite{LP} to obtain a new proof of a classical theorem of Fricke, Klein, and Vogt.

For example, closing off each term in Example \ref{ex:looptwo} gives the recurrence
    \begin{multline}\label{eq:fourtermcf2}
    \tikz[scale=1.5,yscale=.7]{
        \draw[trivalent]
            (0,0)to[bend left=80](0,1)node[leftlabel,pos=.75]{$a$}
            (0,0)to[bend right=80](0,1)node[rightlabel,pos=.75]{$b$};
        \draw[trivalent]
            (0,0)to[bend right](.5,-.3)
            (.5,-.3)to[bend right=80](.5,1.3)node[rightlabel,pos=.7]{$c$}
            (.5,1.3)to[bend right=50](0,1);
        \draw[blue,trivalent]
            (0,.15)to[bend left=80](0,.85)to[bend left=80](0,.15)node[leftlabel,pos=.5]{$\gamma$};
    }
    =
    \tikz[scale=1.5,yscale=.7]{
        \draw[trivalent]
            (0,0)to[bend left=80](0,1)node[leftlabel,pos=.75,scale=.75]{$a+1$}
            (0,0)to[bend right=80](0,1)node[rightlabel,pos=.75,scale=.75]{$b+1$};
        \draw[trivalent]
            (0,0)to[bend right](.5,-.3)
            (.5,-.3)to[bend right=80](.5,1.3)node[rightlabel,pos=.7,scale=.75]{$c$}
            (.5,1.3)to[bend right=50](0,1);
    }
    +\tfrac{\mfe_a(b,c)^2}{b(b+1)}
    \tikz[scale=1.5,yscale=.7]{
        \draw[trivalent]
            (0,0)to[bend left=80](0,1)node[leftlabel,pos=.75,scale=.75]{$a+1$}
            (0,0)to[bend right=80](0,1)node[rightlabel,pos=.75,scale=.75]{$b-1$};
        \draw[trivalent]
            (0,0)to[bend right](.5,-.3)
            (.5,-.3)to[bend right=80](.5,1.3)node[rightlabel,pos=.7,scale=.75]{$c$}
            (.5,1.3)to[bend right=50](0,1);
    }\\
    +\tfrac{\mfe_b(a,c)^2}{a(a+1)}
    \tikz[scale=1.5,yscale=.7]{
        \draw[trivalent]
            (0,0)to[bend left=80](0,1)node[leftlabel,pos=.75,scale=.75]{$a-1$}
            (0,0)to[bend right=80](0,1)node[rightlabel,pos=.75,scale=.75]{$b+1$};
        \draw[trivalent]
            (0,0)to[bend right](.5,-.3)
            (.5,-.3)to[bend right=80](.5,1.3)node[rightlabel,pos=.7,scale=.75]{$c$}
            (.5,1.3)to[bend right=50](0,1);
    }
    +\tfrac{\mfe_c(a,b)^2(\mfe(a,b,c)+1)^2}{a(a+1)b(b+1)}
    \tikz[scale=1.5,yscale=.7]{
        \draw[trivalent]
            (0,0)to[bend left=80](0,1)node[leftlabel,pos=.75,scale=.75]{$a-1$}
            (0,0)to[bend right=80](0,1)node[rightlabel,pos=.75,scale=.75]{$b-1$};
        \draw[trivalent]
            (0,0)to[bend right](.5,-.3)
            (.5,-.3)to[bend right=80](.5,1.3)node[rightlabel,pos=.7,scale=.75]{$c$}
            (.5,1.3)to[bend right=50](0,1);
    }.
    \end{multline}
Note the similarity to the formula in Theorem \ref{t:ranktworecurrencex}.

\subsubsection{Application to Rank Three Central Functions}

The left-associative rank three central functions are
    $$
    \chi_{a,b,c,d,e,f}=
    \tikz[scale=1.4]{
        \draw[trivalent]
            (0,0)to[bend left=80]node[small matrix]{$\xb_1$}(0,1)node[leftlabel,pos=.8]{$a$}
            (0,0)to[bend right=80]node[small matrix]{$\xb_2$}(0,1)node[rightlabel,pos=.8]{$b$}
            (0,0)to[bend right=20](.5,-.2)node[bottomlabel,pos=.5]{$e$}
            to[bend right=80]node[small matrix]{$\xb_3$}(.5,1.2)node[rightlabel,pos=.75]{$c$}
            to[bend right=20](0,1)node[toplabel,pos=.5]{$f$}
            (.5,-.2)to[bend right=20](1,-.4)
            to[bend right=80](1,1.4)node[rightlabel,pos=.75]{$d$}
            to[bend right=20](.5,1.2);
    }
    $$
There are six simple loops in the diagram:
    \begin{align*}
        (a,b) \leftrightarrow \Tr{\xb_1\xb_2^{-1}}, &\qquad (c,d) \leftrightarrow \Tr{\xb_3}, \\
        (a,e,d,f) \leftrightarrow \Tr{\xb_1}, &\qquad (b,e,d,f) \leftrightarrow \Tr{\xb_2}, \\
        (a,e,c,f) \leftrightarrow \Tr{\xb_1\xb_3^{-1}}, &\qquad (b,e,c,f) \leftrightarrow \Tr{\xb_2\xb_3^{-1}}.
    \end{align*}
However, as indicated in the appendix, up to seven variables may be required in the expansion of rank three central functions, so the simple terms do not suffice to compute all central functions. This case will be treated in detail in section \ref{ss:rank3combinatorial}.

\ifnum0=1
    $$\chi_{a_i,c_i}=x_i\cdot\chi_{a_i-1,c_i-1}
            -(-1)^{\beta_i^1+\beta_i^2}\left(\frac{\alpha_i^1\alpha_i^2}{c_i(c_i+1)}\chi_{a_i-2,c_i}
                +\frac{\gamma_i^1\gamma_i^2}{a_i(a_i+1)}\chi_{a_i,c_i-2}\right)
            -\frac{\beta_i^1\beta_i^2(\delta_i^1+1)(\delta_i^2+1)}{a_i(a_i+1)c_i(c_i+1)}\chi_{a_i-2,c_i-2}.$$
\fi

\subsubsection{Application to General Central Functions}
In section \ref{s:centralfunctions}, the left associative central functions were shown to be given by
    $$
    \chi_{\vi,\vec m,\vec p} \equiv
    \tikz[scale=1.2]{
        \draw[trivalent]
            (0,0)to[bend left=80]node[small matrix]{$X_1$}(0,1)node[leftlabel,pos=.8]{$i_1$}
            (0,0)to[bend right=80]node[small matrix]{$X_2$}(0,1)node[rightlabel,pos=.8]{$i_2$}
            (0,0)to[bend right=20](.5,-.2)node[bottomlabel,pos=.5]{$m_1$}
            to[bend right=80]node[small matrix]{$X_3$}(.5,1.2)node[rightlabel,pos=.8]{$i_3$}
            to[bend right=20](0,1)node[toplabel,pos=.5]{$p_1$}
            (.5,-.2)to[bend right=20](1,-.4)node[bottomlabel,pos=.5]{$m_2$}
            (1,1.4)to[bend right=20](.5,1.2)node[toplabel,pos=.5]{$p_2$};
        \draw[draw=none](1.25,0)--(2.25,-.2)node[pos=.2]{.}node[pos=.5]{.}node[pos=.8]{.};
        \draw[draw=none](1.25,1)--(2.25,1.2)node[pos=.2]{.}node[pos=.5]{.}node[pos=.8]{.};
        \draw[trivalent,shift={(.5,0)}]
            (1.5,-.6)to[bend right=20](2,-.8)node[bottomlabel,pos=.4]{$m_{r-2}$}
            to[bend right=80]node[small matrix]{$X_r$}(2,1.8)node[leftlabel,pos=.8]{$i_r$}
            to[bend right=20](1.5,1.6)node[toplabel,pos=.6]{$p_{r-2}$}
            (2,-.8)to[bend right=20](2.5,-1)
            to[bend right=80](2.5,2)node[rightlabel,pos=.4]{$m_{r-1}$}
            to[bend right=20](2,1.8);
    }.
    $$

\begin{theorem}\label{simplelooptheorem}
    There are $\frac{r(r+1)}{2}$ possible simple loops in the general rank $r$ left-associative central function. Multiplication by any of these possible loops gives rise to a recurrence as described in Theorem \ref{t:simplerecurrence}.
\begin{proof}
    Simple loops must pass through precisely two of the edges labeled by $i_l$, implying that there are $\binom{r+1}{2}=\frac{r(r+1)}{2}$ possible simple loops.
\end{proof}
\end{theorem}
The simple recurrences given in this way correspond to the following subset of generators of $\C[\R_r]^G$:  $\{\Tr{\xb_1},...,\Tr{\xb_r}\}\cup\{\Tr{\xb_i\xb_j}\ |\ 1\leq i<j \leq r\}$ (see appendix for a full description of minimal generators).

%% file: nonsimplerecurrence.tex
\subsection{Non-Simple Loop Recurrences}\label{s:nonsimple}

\begin{figure}[h]
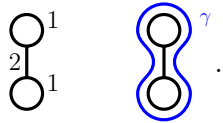

    \tikz[trivalent,every node/.style={basiclabel}]{
        \draw(0,-.1)circle(.3)(0,1.1)circle(.3)(0,.2)to node[auto]{$2$}(0,.8);
        \node at(.5,.1){$1$};\node at(.5,1.3){$1$};}
    \qquad
    \tikz[trivalent,every node/.style={basiclabel}]{
        \draw(0,-.1)circle(.3)(0,1.1)circle(.3)(0,.2)to(0,.8);
        \foreach\xa in{1,-1}{
            \draw[blue,xscale=\xa](0,-.6)to[out=180,in=-90](-.5,-.1)to[wavyup](-.2,.5)to[wavyup](-.5,1.1)to[out=90,in=180](0,1.6);
        }
        \node[blue,rightlabel]at(.6,1.3){$\gamma$};
    }.
    \caption{Trace diagram with no simple recurrence (left). The non-simple loop which can be used to produce a recurrence (right).}\label{fig:barbellcase}
    \end{figure}

As mentioned in the previous section, multiplication by simple loops does not provide sufficient recurrence relations to reduce any trace diagram to its simplest pieces. In the rank three case, for example, multiplication by $\Tr{\xb_1\xb_2^{-1}\xb_3}$ requires a non-simple loop:
    $$
    \tikz[scale=1.4,xscale=1.8]{
        \draw[trivalent]
            (0,0)to[bend left=50]node[small matrix]{$\xb_1$}(0,1)to[bend left=50]node[small matrix]{$\xb_2$}(0,0)
            (0,0)to[bend right=20](.5,-.2)to[bend right=80]node[small matrix]{$\xb_3$}(.5,1.2)to[bend right=20](0,1)
            (.5,-.2)to[bend right=20](1,-.4)to[bend right=80](1,1.4)to[bend right=20](.5,1.2);
        \draw[blue,trivalent]
            (-.1,-.1)to[bend left=80]node[small matrix]{$\xb_1$}(0,1.1)to[bend left=80]node[small matrix,pos=.55]{$\xb_2$}(.2,0)
            to[bend right=20](.5,-.1)to[bend right=60]node[small matrix]{$\xb_3$}(.4,1.1)to[bend left=40](1.1,1.5)
            to[bend left=80](1.1,-.5)to[bend left=30](.5,-.3)to[bend left=20](-.1,-.1);
    }
    $$
While Corollary \ref{c:simplerecurrence} can be applied to any simple loop in a trace diagram, it is not sufficient to compute the value of an arbitrary diagram. The simplest case where it fails is the ``barbell'' depicted in Figure \ref{fig:barbellcase}. There is no simple loop recurrence, since subtracting one from either loop produces a non-admissible diagram. However, the non-simple loop depicted on the right of Figure \ref{fig:barbellcase} does provide a recurrence.

Define
    $$
    \tilde\chi_{a,c}^b({\bf X},{\bf Y})=
    \tikz[trivalent]{
        \draw(0,-.3)circle(.4)(0,1.3)circle(.4)(0,.1)to(0,.9);
        \node[basiclabel]at(.5,.1){$a$};\node[basiclabel]at(.5,1.6){$c$};
        \node[rightlabel]at(0,.5){$b$};
        \node[small matrix]at(-.4,-.3){$\bf X$};\node[small matrix]at(-.4,1.3){$\bf Y$};}
    $$
When applied to the simple loops in the diagram, \eqref{eq:loopone} implies
    \begin{align*}
        \Tr{\bf X} \cdot \tilde\chi_{a,c}^b &=
            \tilde\chi_{a+1,c}^b - \nFus ba{a-1}a{a-1} \cdot \tilde\chi_{a-1,c}^b \\
        \Tr{\bf Y} \cdot \tilde\chi_{a,c}^b &=
            \tilde\chi_{a,c+1}^b - \nFus bc{c-1}c{c-1} \cdot \tilde\chi_{a,c-1}^b.
    \end{align*}
These equations can be rearranged to provide recurrences for the $\bf X$ and $\bf Y$ loops. The third recurrence corresponds to the multiplication by $\Tr{\bf XY}$:

\begin{proposition}\label{prop:barbellrecurrence}
    $$
    \tikz[trivalent]{
        \foreach\xa in{1,-1}{
            \draw[blue,xscale=\xa](0,-.9)to[out=180,in=-90](-.6,-.3)to[wavyup](-.2,.5)to[wavyup](-.6,1.3)to[out=90,in=180](0,1.9);}
        \draw(0,-.3)circle(.4)(0,1.3)circle(.4)(0,.1)to(0,.9);
        \node[basiclabel]at(.2,-.3){$a$};\node[basiclabel]at(.2,1.3){$c$};
        \node[coordinate,pin={[basiclabel,pin edge={darkgray}]right:$b$}]at(0,.8){};
        \node[small matrix]at(-.4,-.3){$\bf X$};\node[small matrix]at(-.4,1.3){$\bf Y$};
        \node[blue,rightlabel]at(.5,1.6){$1$};
    }
    =\! \!\!\! \!\!\! \!\!\! \!\!\!\!\! \!\!\!\!
    \sum_{\substack{a'\in\iadm{a,1}\\
            b'\in\iadm{b,1}\cap\iadm{a,a'}\\
            c'\in\iadm{c,1}\cap\iadm{c,b'}\\
            b''\in\iadm{b',1}\cap\iadm{a,a'}\cap\iadm{c,c'}}}\hspace{-1cm}
    \mfs_{a'}(1,a)\mfs_{c'}(1,c)\nFus aa{a'}b{b'}\nFus cc{c'}b{b'}\nFus{a'}a{a'}{b'}{b''}\nFus{c'}c{c'}{b'}{b''}
    \tikz[trivalent]{
        \draw(0,-.3)circle(.4)(0,1.3)circle(.4)(0,.1)to(0,.9);
        \node[basiclabel]at(.5,.1){$a'$};\node[basiclabel]at(.5,1.6){$c'$};
        \node[rightlabel]at(0,.5){$b''$};
        \node[small matrix]at(-.4,-.3){$\bf X$};\node[small matrix]at(-.4,1.3){$\bf Y$};}.
    $$
\begin{proof}
    The computation is a straightforward application of Lemma \ref{l:vertexglue} and the bubble identity \eqref{eq:bubbleidentity}. For simplicity, the matrices are not shown on the diagram; nothing changes if they are included.
    \begin{align*}
    \tikz[trivalent]{
        \foreach\xa in{1,-1}{
            \draw[blue,xscale=\xa](0,-.9)to[out=180,in=-90](-.6,-.3)to[wavyup](-.2,.5)to[wavyup](-.6,1.3)to[out=90,in=180](0,1.9);}
        \draw(0,-.3)circle(.4)(0,1.3)circle(.4)(0,.1)to(0,.9);
        \node[basiclabel]at(.15,-.25){$a$};\node[basiclabel]at(.15,1.35){$c$};
        \node[coordinate,pin={[basiclabel,pin edge={darkgray}]right:$b$}]at(0,.8){};
        \node[blue,rightlabel]at(.5,1.6){$1$};
    }
    =
    \hspace{-.5cm}
    \sum_{\substack{a'\in\iadm{a,1}\\
            b'\in\iadm{b,1}\cap\iadm{a,a'}\\
            c'\in\iadm{c,1}\cap\iadm{c,b'}}}
    \hspace{-.5cm}
    &    \mfs_{b'}(1,b)\Fus aa{a'}b{b'}\mfs_{b'}(1,b)\Fus cc{c'}b{b'}\mfb_{b'}(1,b)
    \tikz[trivalent]{
        \draw[blue](0,-.9)to[out=0,in=-90](.6,-.3)to[wavyup](.2,.5)to[wavyup](.6,1.3)to[out=90,in=0](0,1.9)
            (0,-.9)to[out=180,in=-90](-.6,-.5)to[out=90,in=-135](-.4,-.3)
            (0,1.9)to[out=180,in=90](-.6,1.5)to[out=-90,in=135](-.4,1.3);
        \draw(0,-.3)circle(.4)(0,1.3)circle(.4)(0,.1)to(0,.9);
        \node[basiclabel]at(.15,-.25){$a$};\node[basiclabel]at(.15,1.35){$c$};
        \node[leftlabel]at(-.2,.15){$a'$};\node[leftlabel]at(-.3,.9){$c'$};
        \node[coordinate,pin={[basiclabel,pin edge={darkgray}]right:$b'$}]at(0,.8){};
        \node[blue,rightlabel]at(.5,1.6){$1$};
    } \\
    =
    \hspace{-1cm}
    \sum_{\substack{a'\in\iadm{a,1}\\
            b'\in\iadm{b,1}\cap\iadm{a,a'}\\
            c'\in\iadm{c,1}\cap\iadm{c,b'}\\
            b''\in\iadm{b',1}\cap\iadm{a,a'}\cap\iadm{c,c'}}}
    \hspace{-1cm}
    & \Fus aa{a'}b{b'}\Fus cc{c'}b{b'}\mfb_{b'}(1,b)
        \mfs_{b''}(1,b')\Fus {a'}a{a'}{b'}{b''}\mfs_{b''}(1,b')\Fus{c'}c{c'}{b'}{b''}\mfb_{b''}(1,b')
    \tikz[trivalent]{
        \draw[blue]
            (0,-.9)to[out=180,in=-90](-.6,-.5)to[out=90,in=-135](-.4,-.3)
            (0,-.9)to[out=0,in=-90](.6,-.5)to[out=90,in=-45](.4,-.3)
            (0,1.9)to[out=180,in=90](-.6,1.5)to[out=-90,in=135](-.4,1.3)
            (0,1.9)to[out=0,in=90](.6,1.5)to[out=-90,in=45](.4,1.3);
        \draw(0,-.3)circle(.4)(0,1.3)circle(.4)(0,.1)to(0,.9);
        \node[toplabel]at(0,-.7){$a$};\node[bottomlabel]at(0,1.7){$c$};
        \node[leftlabel]at(-.2,.15){$a'$};\node[leftlabel]at(-.3,.9){$c'$};
        \node[rightlabel]at(.3,.1){$a'$};\node[rightlabel]at(.3,1.1){$c'$};
        \node[rightlabel]at(0,.5){$b''$};
        \node[blue,rightlabel]at(.6,-1){$1$};
        \node[blue,rightlabel]at(.5,1.6){$1$};
    } \\
    =
    \hspace{-1cm}
    \sum_{\substack{a'\in\iadm{a,1}\\
            b'\in\iadm{b,1}\cap\iadm{a,a'}\\
            c'\in\iadm{c,1}\cap\iadm{c,b'}\\
            b''\in\iadm{b',1}\cap\iadm{a,a'}\cap\iadm{c,c'}}}
    \hspace{-1cm}
    &    \mfs_{a'}(1,a)\mfs_{c'}(1,c)\nFus aa{a'}b{b'}\nFus cc{c'}b{b'}\nFus{a'}a{a'}{b'}{b''}\nFus{c'}c{c'}{b'}{b''}
    \tikz[trivalent,every node/.style={basiclabel}]{
        \draw(0,-.1)circle(.3)(0,1.1)circle(.3)(0,.2)to node[auto]{$b''$}(0,.8);
        \node at(.5,.1){$a'$};\node at(.5,1.3){$c'$};
    }
    \end{align*}
\end{proof}
\end{proposition}

The explicit computation of these functions is shown in Table \ref{t:barbell} with $x=\mathrm{tr}(\xb)$, $y=\mathrm{tr}(\mathbf{Y})$, and $z=\mathrm{tr}(\xb\mathbf{Y})$. We omit the cases where $b=0$, since they can be obtained as a product of two rank one central functions.

{\renewcommand\arraystretch{1.5}
  \begin{table}[ht]
    \begin{center}
      \begin{tabular}{|c|c|c|}
        \hline
        {\bf a+b+c} & $\tilde\chi_{a,c}^b$ & formula\\
        \hline\hline
        $\bf 4$ & $\tilde\chi_{1,1}^2$ & $z-\frac{x y}{2}$\\
        \hline
        $\bf 5$ & $\tilde\chi_{2,1}^2$ & $x z-\frac{x^2 y}{2}$\\
                & $\tilde\chi_{1,2}^2$ & $y z-\frac{x y^2}{2}$\\
        \hline
        $\bf 6$ & $\tilde\chi_{3,1}^2$ & $x^2 z-\frac{x^3 y}{2} + \frac{x y}{3}-\frac{2 z}{3}$ \\
                & $\tilde\chi_{2,2}^2$ & $x y z-\frac{1}{2} x^2 y^2$ \\
                & $\tilde\chi_{1,3}^2$ & $y^2 z-\frac{x y^3}{2} + \frac{x y}{3}-\frac{2 z}{3}$ \\
        \hline
        $\bf 7$ & \emph{4 cases} & $\cdots$ \\
        \hline
        $\bf 8$ & $\tilde\chi_{2,2}^4$ & $z^2-x y z+\frac{x^2 y^2}{6}+\frac{x^2}{3}+\frac{y^2}{3}-\frac{4}{3}$ \\
                & \emph{5 more cases} & $\cdots$ \\
        \hline
      \end{tabular}
    \medskip
    \caption{Barbell Functions.}\label{t:barbell}
    \end{center}
  \end{table}
}

\subsubsection{The General Case}

Other non-simple recurrence formulas may be similarly derived. In every case, the coefficients will be signed summations over the normalized fusion coefficients, with the sign depending on the particular configuration of the graph.

%% file: rank3.tex
\subsection{Tensorial Algorithm}
We now heuristically describe an effective algorithm to compute central functions for ranks 1, 2, and 3 using a ``tensorial contraction'' method.  This algorithm has been successfully implemented in {\it Mathematica}, and a semi-documented ``notebook'' is available to the reader (at the first author's website \cite{Lweb}).

Let $t:\R_r\to \C^{N_r}$ be the mapping given by specifying $N_r$ minimal generators for $\C[\R_r]^G$, let $\X_r=t(\R_r)$ be the image of $t$.  The projection $\R_r\to \X_r$ is dual to the inclusion $\C[\R_r]^G\subset \C[\R_r]$. Denote by $\overline{[\rho]}$ the equivalence class of $\rho\in \R_r$ where $\rho\sim \rho'$ if and only if $\overline{G\rho}\cap\overline{G\rho'}\not=\emptyset$, that is their conjugation orbit closures intersect non-trivially. $\X_r$ can be understood as the space of such equivalence classes.  In the case $r=3$ this mapping is given by
$\overline{[\rho]}\mapsto$ $$\left(\Tr{\rho(\xt_1)},\Tr{\rho(\xt_2)},\Tr{\rho(\xt_3)},\Tr{\rho(\xt_1\xt_2)},\Tr{\rho(\xt_1\xt_3)},\Tr{\rho(\xt_2\xt_3)},\Tr{\rho(\xt_1\xt_2\xt_3)}\right).$$
In \cite{G9} an explicit global slice to this projection is constructed.  Using this slice and an iterative application of the work in \cite{LP} we compute rank 3 central functions tensorially and describe this algorithm below.

To ease the notation let $t_k=\Tr{\xb_k}=\Tr{\rho(\xt_k)}$, $t_{ij}=\Tr{\xb_i\xb_j}=\Tr{\rho(\xt_i\xt_j)}$, and $t_{ijk}=\Tr{\xb_i\xb_j\xb_k}=\Tr{\rho(\xt_i\xt_j\xt_k)}$.

The projection $\pi:\X_3\to \C^6$ given by $\overline{[\rho]}\mapsto (t_1,t_2,t_3,t_{12},t_{13},t_{23})$ is a branched double cover.  In particular, there exist polynomials in $P,Q\in\C[t_1,t_2,t_3,t_{12},t_{13},t_{23}]$ so that $$\C[\X_3]\approx \C[t_1,t_2,t_3,t_{12},t_{13},t_{23}][t_{123}]/(t_{123}^2-Pt_{123}+Q).$$  The two roots of the irreducible generator of the ideal are related by the formula:
$$t_{132}= -t_1 t_2 t_3+t_{12} t_3+t_2
   t_{13}+t_1 t_{23}-t_{123}.$$

This provides the formula for $P$; the formula for $Q$ is
$$t_{123} t_{132} = (t_{1}^2
+ t_{2}^2+ t_{3}^2) + (t_{12}^2 + t_{23}^2 + t_{13}^2) -
(t_{1}t_{2}t_{12} + t_{2}t_{3}t_{23} + t_{3}t_{1}t_{13}) + t_{12}t_{23}t_{13} - 4.$$

To construct a slice we must construct a triple of matrices $(\xb_1,\xb_2,\xb_3)$ for every 7-tuple in the image of $t$.

Let $\xb_1=\left(
\begin{array}{ll}
 t_1 & -1 \\
 1 & 0
\end{array}
\right)$ and $\xb_2=\left(
\begin{array}{ll}
 0 & w \\
 -\frac{1}{w} & t_2
\end{array}
\right)$, where $w+1/w=t_{12}$.  Then letting $$\xb_3(s)=\left(
\begin{array}{ll}
 s \left(\frac{1}{w}-w\right)+t_3 & s \left(w
   t_1-t_2\right)+\frac{w \left(w
   \left(t_{13}-t_1
   t_3\right)+t_{23}\right)}{w^2-1} \\
 s \left(\frac{t_1}{w}-t_2\right)+\frac{-t_1
   t_3+t_{13}+w t_{23}}{w^2-1} & s
   \left(w-\frac{1}{w}\right)
\end{array}
\right),$$ where $s$ is a solution to the equation $\mathrm{det}(\xb_3(s))-1=0$, will give the desired slice (see \cite{G9}).

\subsubsection{Step 1}
Define tensor products and duality scalar products.  The scalar product satisfies

$${\sf n}^*_{n-k}({\sf n}_{n-l})=\frac{(n-k)!k!}{n!}\delta_{kl}=\raisebox{2pt}{$\delta_{kl}$}\!\Big/\!\raisebox{-2pt}{$\tbinom{n}{k}$},$$ for ${\sf n}_{n-k}\in V_n$.

\subsubsection{Step 2}
Define an algorithm which determines the form of a linear representation of an element in $\SL$ on the symmetric power $V_n$.

This comes from
\begin{align*}
{\sf n}^*_{n-k}(\xb\cdot {\sf n}_{n-l})
&=  {\sf n}^*_{n-k}\left((x_{11} e_1+x_{21} e_2)^{n-l}(x_{12} e_1 + x_{22} e_2)^l\right)\\
&=\sum_{\substack{i+j=k\\0 \le i \le n-l \\ 0 \le j \le l}}
\tbinom{n}{k}^{-1}\tbinom{n-l}{i}\tbinom{l}{j}x_{11}^{n-l-i}x_{12}^{l-j}x_{21}^ix_{22}^j.
\end{align*}

For instance, if $g=\left(
\begin{array}{ll}
 a & b \\
 c & d
\end{array}
\right)$, then the induced action on $V_2$ is given by $\left(
\begin{array}{lll}
 a^2 & a b & b^2 \\
 a c & \frac{1}{2} (b c+a d) & b d \\
 c^2 & c d & d^2
\end{array}
\right)$, and on $V_3$ by

$$\left(
\begin{array}{llll}
 a^3 & a^2 b & a b^2 & b^3 \\
 a^2 c & \frac{1}{3} \left(d a^2+2 b c a\right)
   & \frac{1}{3} \left(c b^2+2 a d b\right) &
   b^2 d \\
 a c^2 & \frac{1}{3} \left(b c^2+2 a d c\right)
   & \frac{1}{3} \left(a d^2+2 b c d\right) & b
   d^2 \\
 c^3 & c^2 d & c d^2 & d^3
\end{array}
\right).$$

Upon contracting tensors, we will need to know the polynomial matrix coefficients of the action on a given $V_n$.  This routine will allow us to read off such matrix coefficients, since the pairing for ${\sf n}^*_{n-k}(\xb\cdot {\sf n}_{n-l})$ just becomes the $n-k+1,n-l+1$ entry of the symmetrized matrix.

\subsubsection{Step 3}
Given a triple $(a,b,c)$ determine all admissible 6-tuples, and mark them by an enumeration of the multiplicities.  We iteratively decompose tensor triples $V_a\otimes V_b\otimes V_c$ using a left associative iteration algorithm.  In other words, we decompose $(V_a\otimes V_b)\otimes V_c$, using the decomposition of $V_a\otimes V_b\approx \sum V_{a+b-2k}$.  Then for each allowed value of $0\leq k\leq a+b$, we decompose $V_{a+b-2k}\otimes V_c$ in the same fashion.  We know that $(a,b,a+b-2k)$ is admissible if and only if  $|a-b| \leq  a+b-2k \leq  a+b$.

For instance $V_d$ injects into $V_3\otimes V_2\otimes V_2$ if and only if $d\in \{7, 5, 3, 5, 3, 1, 3, 1\}$. The multiplicity arises since

\begin{align*}
(V_3\otimes V_2)\otimes V_2&\approx (V_5\oplus V_3 \oplus V_1)\otimes V_2 \approx V_5\otimes V_2\oplus V_3\otimes V_2 \oplus V_1\otimes V_2\\
&\approx V_7\oplus V_5\oplus V_3 \oplus V_5\oplus V_3\oplus V_1 \oplus V_3 \oplus V_1.
\end{align*}

\subsubsection{Step 4}
Define the injection of basic elements.  This is done using the formula for the mapping $\iota:V_c\to V_a\otimes V_b$:
\begin{align*}
  \tbinom{c}{k}{\sf c}_k
    &\longmapsto\sum_{\substack{0\le i\le\beta\\0\le j\le\alpha\\0\le m\le\gamma\\i+j=k}}%
      \tbinom{\beta}{i}{\sf a}_i
        \otimes\left[(-1)^m\tbinom{\gamma}{m}{\sf a}_{\gamma-m}\otimes{\sf b}_m\right]
        \otimes\tbinom{\alpha}{j}{\sf b}_j\\ %
    &\longmapsto\sum_{\substack{0\le i\le\beta\\0\le j\le\alpha\\0\le m\le\gamma\\i+j=k}}%
      (-1)^m\tbinom{\beta}{i}\tbinom{\alpha}{j}\tbinom{\gamma}{m}{\sf a}_{i+\gamma-m}\otimes{\sf b}_{j+m},%
\end{align*}
and iterating for each summand.

Also since we are using left associations, i.e. grouping $V_a$ and $V_b$ in $V_a\otimes V_b\otimes V_c$, we first inject $V_d\hookrightarrow V_{a+b-2k}\otimes V_c$ and then inject $V_{a+b-2k}\otimes V_c\hookrightarrow V_a\otimes V_b\otimes V_c$ using the injection $V_{a+b-2k}\hookrightarrow V_a\otimes V_b$ tensored with the identity.

Likewise, we have an injection for the dual injection $V_d^*\hookrightarrow V_a^*\otimes V_b^*\otimes V_c^*$.

\subsubsection{Step 5}

The central function begins in $\mathrm{End}(V_d)^{\SL}$ and then is mapped to $(V^*_d \otimes V_d)^{\SL}$ by
$$\sum_{k=0}^d{\sf d}_k({\sf d}_k)^T\mapsto\sum_{k=0}^c\tbinom{d}{k}\bs{\sf d}{k}{k}.$$
The central function is then determined by the composite injections $$\mathrm{End}(V_d)^{\SL}\approx (V^*_d \otimes V_d)^{\SL} \hookrightarrow (V^*_a\otimes V^*_b \otimes V^*_c\otimes V_a\otimes V_b\otimes V_c)^{\SL}.$$  Once we realize the explicit form of the injections $V_d\hookrightarrow V_a\otimes V_b\otimes V_c$ we can write the central function $\sum_{k=0}^d{\sf d}_k({\sf d}_k)^T$ as a central tensor in terms of $V_a\otimes V_b\otimes V_c$.  We note that the coefficients are a bit delicate here:  first one chooses to include $\tbinom{d}{k}$ on only one the vectors ${\sf d_k}$ or ${\sf d}^*_k$, but not both, for each summand of the central function.  Second one must make sure to include a $1/\tbinom{d}{k}$ when including each factor of each summand (both the vector and its dual) into a tensor $V_a\otimes V_b$.  We observe that the rank 2 coefficients cancel, but for rank 3 they generally do not cancel since we are iterating injections.  Lastly, we have included the dual pairing binomial coefficients in our expression of the symmetrization of a generic matrix, so we do not need to include any further binomial coefficients.

\subsubsection{Step 6}

Once we have the central function as a linear combination of tensors, we must ``contract the tensor'' using the mapping that associates a linear combination of tensor (of the type we have been considering) to a polynomial function in the coordinate ring $\C[\R_r]$.  Recall that this is determined by mapping
$(e^*_{i_1}\otimes e^*_{i_2}\otimes \cdots \otimes e^*_{i_r} )\otimes (e_{j_1}\otimes e_{j_2}\otimes \cdots \otimes e_{j_r} )$ to the polynomial function
$$(\xb_1, \xb_2,...,\xb_r) \mapsto e^*_{i_1}(\xb_1\cdot e_{j_1})e^*_{i_2}(\xb_2\cdot e_{j_2})\cdots e^*_{i_r} (\xb_r\cdot e_{j_r} ).$$  Once a generic matrix $\xb_i$ is mapped to an automorphism of $V_d$, call it $\mathrm{Sym}^d(\xb_i)$.  Then $e^*_{i_1}(\xb_1\cdot e_{j_1})$ is just the $(i_1+1, j_1+1)$ entry of $\mathrm{Sym}^d(\xb_i)$.

\subsubsection{Step 7}

Lastly, using the Goldman slice (see \cite{G9}) $$\X_3=\SL^{\times 3}\aq \SL \hookrightarrow \R_3=\SL^{\times 3}$$ to the categorical projection $\R_3\to\X_3$, we express the invariant polynomials (which are polynomials in the generic matrix entries) from the tensorial contraction in terms of the seven invariants: $$\{\Tr{\xb_1},\Tr{\xb_2},\Tr{\xb_3},\Tr{\xb_1\xb_2},\Tr{\xb_1\xb_3},\Tr{\xb_2\xb_3},\Tr{\xb_1\xb_2\xb_3}\}.$$

It is only this last step that does not generalize to $r\geq 4$, since we do not have a slice.  However, for $r\geq 4$, we can still contract the tensor (Steps 1- 6) and then use a Gr\"obner Basis algorithm to express the polynomials in terms of traces (minimally in fact since minimal generators are known for all $r\geq 1$).

\subsection{Combinatorial Algorithm}\label{ss:rank3combinatorial}
Trace diagrams allow a more combinatorial approach to computing the central functions. The key point is that any central function can be reduced in terms of simpler central functions. As discussed in section \ref{simplerecurrence-section}, a recurrence formula exists for each loop in the diagrammatic depiction of the central function.

In this section, we show how this process works for $r=3$ diagrams of the form:
    $$
    \chi_{a,b,c,d,e,f}=
    \tikz[scale=1.4]{
        \draw[trivalent]
            (0,0)to[bend left=80]node[small matrix]{$\xb_1$}(0,1)node[leftlabel,pos=.8]{$a$}
            (0,0)to[bend right=80]node[small matrix]{$\xb_2$}(0,1)node[rightlabel,pos=.8]{$b$}
            (0,0)to[bend right=20](.5,-.2)node[bottomlabel,pos=.5]{$e$}
            to[bend right=80]node[small matrix]{$\xb_3$}(.5,1.2)node[rightlabel,pos=.75]{$c$}
            to[bend right=20](0,1)node[toplabel,pos=.5]{$f$}
            (.5,-.2)to[bend right=20](1,-.4)
            to[bend right=80](1,1.4)node[rightlabel,pos=.75]{$d$}
            to[bend right=20](.5,1.2);
    }
    $$
We call the sum of $a+b+c$ the \emph{fundamental order}.  Then the algorithm reduces the diagram into a sum over diagrams in a lower fundamental order.  The algorithm terminates with the base case $\chi_{0,0,0,0,0,0}=1$. The process uses eight different recurrence formulas, corresponding to the following trace variables:
    \begin{multline*}
    \{\tr(\xb_1), \tr(\xb_2), \tr(\xb_3), \tr(\xb_1 \xb_2^{-1}), \tr(\xb_1 \xb_3^{-1}), \tr(\xb_2 \xb_3^{-1}),\\ \tr(\xb_1 \xb_2^{-1} \xb_3), \tr(\xb_1 \xb_3 \xb_2^{-1})\}.
    \end{multline*}
The final two recurrences, which correspond to nonsimple loops in the diagram, require particular attention.

\subsubsection{Algorithm Overview}
As mentioned in section \ref{simplerecurrence-section}, the recurrence corresponding to a loop with $n$ edges contains a maximum of $2^n$ terms. For that reason, it is more efficient to begin the computation with the shortest loops in the figure. The algorithm implemented in Mathematica begins with the following steps:
\begin{enumerate}
    \item If $\mathfrak{e}_e(a,b)>0$ and $\mathfrak{e}_f(a,b)>0$, reduce along the $(a,b)$ loop, corresponding to $\tr(\xb_1 \xb_2^{-1})$;
    \item otherwise, if $\mathfrak{e}_e(c,d)>0$ and $\mathfrak{e}_f(c,d)>0$, reduce along the $(c,d)$ loop, corresponding to $\tr(\xb_3)$;
    \item otherwise, if $\mathfrak e_b(a,e)>0$, $\mathfrak e_b(a,f)>0$, $\mathfrak e_c(d,e)>0$, and $\mathfrak e_c(d,f)>0$, reduce along the $(a,e,d,f)$ loop, corresponding to $\tr(\xb_1)$;
    \item otherwise, if $\mathfrak e_a(b,e)>0$, $\mathfrak e_a(b,f)>0$, $\mathfrak e_c(d,e)>0$, and $\mathfrak e_c(d,f)>0$, reduce along the $(b,e,d,f)$ loop, corresponding to $\tr(\xb_2)$;
    \item otherwise, if $\mathfrak e_a(b,e)>0$, $\mathfrak e_a(b,f)>0$, $\mathfrak e_d(c,e)>0$, and $\mathfrak e_d(c,f)>0$, reduce along the $(b,e,c,f)$ loop, corresponding to $\tr(\xb_2 \xb_3^{-1})$;
    \item otherwise, if $\mathfrak e_b(a,e)>0$, $\mathfrak e_b(a,f)>0$, $\mathfrak e_d(c,e)>0$, and $\mathfrak e_d(c,f)>0$, reduce along the $(a,e,c,f)$ loop, corresponding to $\tr(\xb_1 \xb_3^{-1})$;
    \item otherwise, if $e=0$, reduce along the nonsimple $(a,b,e,c,d,e)$ loop, corresponding to $\tr(\xb_1 \xb_3 \xb_2^{-1})$;
    \item otherwise, if $f=0$, reduce along the nonsimple $(a,b,f,c,d,f)$ loop, corresponding to $\tr(\xb_1 \xb_2^{-1} \xb_3)$.
\end{enumerate}

\begin{proposition}
    Any admissible $\chi_{a,b,c,d,e,f}$ can be reduced via one of the above cases.
\begin{proof}
    Consider the following set of edge types:
    \begin{equation}\label{eq:testedges}
        \{\mathfrak e_a(b,e), \mathfrak e_b(a,e), \mathfrak e_a(b,f), \mathfrak e_b(a,f)\}.
    \end{equation}
    If any two of these is zero, then the central function is reducible. To see this, suppose without loss of generality that $\mathfrak e_a(b,e)=0$. If $\mathfrak e_b(a,e)=0$ then $e=\mathfrak e_a(b,e)+\mathfrak e_b(a,e)=0$. If $\mathfrak e_b(a,f)=0$, then $a=b+e$ and $b=a+f$ so that $e=0$ and $f=0$. Finally, if $\mathfrak e_a(b,f)=0$ then $a=b+e=b+f$ and so $\mathfrak e_e(a,b)=2b=\mathfrak e_f(a,b)$. So either $b>0$ and reduction by the $(a,b)$ loop is possible, or $b=0$. In this case, the diagram is reducible by either $\tr(\xb_1)$, $\tr(\xb_1 \xb_3^{-1})$, or $\tr(\xb_3)$.

    The same logic can be applied to the set of edges
    \begin{equation}\label{eq:testedges2}
        \{\mathfrak e_c(d,e), \mathfrak e_d(c,e), \mathfrak e_c(d,f), \mathfrak e_d(c,f)\}.
    \end{equation}
    But note that if at most one of \eqref{eq:testedges} and at most one of \eqref{eq:testedges2} is zero, then one of the four reductions corresponding to a loop of length four is possible.
\end{proof}
\end{proposition}

\subsubsection{Reduction along Simple Loops}
The simple loop reductions may have up to 4 terms in the $\tr(\xb_3)$ and $\tr(\xb_1 \xb_2^{-1})$ cases, and up to 16 terms in the other four cases. While this is a large number of terms, writing out the recurrence is a straightforward process. Each is an immediate consequence of Theorem \ref{t:simplerecurrence}. For example:
\begin{multline*}
    \tr(\xb_1\xb_2^{-1})\chi_{a,b,c,d,e,f} =
         \chi_{a+1,b+1,c,d,e,f}\\
         + \mathfrak s_e(a,b)\mathfrak s_f(a,b)\frac{\mathfrak e_a(b,e) \mathfrak e_a(b,f)}{b(b+1)} \chi_{a+1,b-1,c,d,e,f}\\
         + \mathfrak s_e(a,b)\mathfrak s_f(a,b)\frac{\mathfrak e_b(a,e)\mathfrak e_b(a,f)}{a(a+1)} \chi_{a-1,b+1,c,d,e,f}\\
         + \frac{\mathfrak e_e(a,b)\mathfrak e_f(a,b)(\mathfrak e(a,b,e)+1)(\mathfrak e(a,b,f)+1)}{a(a+1)b(b+1)} \chi_{a-1,b-1,c,d,e,f}.
\end{multline*}
The remainder of the recurrences may be found in a Mathematica notebook written by the second author (available on his website \cite{Pweb}).

\subsubsection{Reduction along Non-Simple Loops}
In the final cases, either $e=0$ or $f=0$, and the central function is topologically equivalent to a barbell. Proposition \ref{prop:barbellrecurrence} could be adapted to this case, but we prefer to show directly how the rank three central functions are related to the barbell functions $\tilde\chi_{a,c}^b$. When $f=0$, the diagram is
    $$
    \chi_{a,b,c,d,e,0} =
    \tikz[scale=1.4]{
        \draw[trivalent]
            (0,0)to[bend left=80]node[small matrix]{$\xb_1$}(0,1)node[leftlabel,pos=.8]{$a$}
            (0,0)to[bend right=80]node[small matrix]{$\xb_2$}(0,1)node[rightlabel,pos=.8]{$b$}
            (0,0)to[bend right=20](.5,-.2)node[bottomlabel,pos=.5]{$e$}
            to[bend right=80]node[small matrix]{$\xb_3$}(.5,1.2)node[rightlabel,pos=.75]{$c$}
            (.5,-.2)to[bend right=20](1,-.4)
            to[bend right=80](1,1.4)node[rightlabel,pos=.75]{$d$}
            to[bend right=20](.5,1.2);
    }
    = (-1)^{\frac e2}\hspace{-12pt}
    \tikz[trivalent]{
        \draw(0,-.5)circle(.6)(0,1.5)circle(.6)(0,.1)to(0,.9);
        \node[rightlabel]at(.5,.1){$c$};\node[rightlabel]at(.6,1.6){$a$};\node[rightlabel]at(0,.5){$e$};
        \node[small matrix]at(-.6,-.5){$\xb_3$};\node[small matrix,ellipse]at(-.6,1.5){$\xb_2^{-1}\xb_1$};}
    = (-1)^{\frac e2} \tilde{\chi}_{c,a}^e(\xb_3,\xb_2^{-1}\xb_1).
    $$
The sign $(-1)^{\frac e2}$ arises from Proposition \ref{p:spinnetsignstrong} since there are $\frac e2$ kinks in the left diagram and none on the right. So $\chi(a,b,c,d,e,0)$ may be obtained from Table \ref{t:barbell} using the substitutions
    \begin{align*}
        x &\mapsto \mathrm{tr}(\xb_3) = t_3,\\
        y &\mapsto \mathrm{tr}(\xb_2^{-1}\xb_1) = t_1t_2-t_{12},\\
        z &\mapsto \mathrm{tr}(\xb_3\xb_2^{-1}\xb_1) = t_{123}+t_1t_2t_3-t_{12}t_3-t_1t_{23}.
    \end{align*}

Similarly, when $e=0$, the diagram is
    $$
    \chi_{a,b,c,d,0,f} =
    \tikz[scale=1.4]{
        \draw[trivalent]
            (0,0)to[bend left=80]node[small matrix]{$\xb_1$}(0,1)node[leftlabel,pos=.8]{$a$}
            (0,0)to[bend right=80]node[small matrix]{$\xb_2$}(0,1)node[rightlabel,pos=.8]{$b$}
            (.5,-.2)to[bend right=80]node[small matrix]{$\xb_3$}(.5,1.2)node[rightlabel,pos=.75]{$c$}
            to[bend right=20](0,1)node[toplabel,pos=.5]{$f$}
            (.5,-.2)to[bend right=20](1,-.4)
            to[bend right=80](1,1.4)node[rightlabel,pos=.75]{$d$}
            to[bend right=20](.5,1.2);
    }
    = (-1)^{\frac f2}\hspace{-12pt}
    \tikz[trivalent]{
        \draw(0,-.5)circle(.6)(0,1.5)circle(.6)(0,.1)to(0,.9);
        \node[rightlabel]at(.5,.1){$a$};\node[rightlabel]at(.6,1.6){$c$};\node[rightlabel]at(0,.5){$f$};
        \node[small matrix,ellipse]at(-.6,-.5){$\xb_1\xb_2^{-1}$};\node[small matrix]at(-.6,1.5){$\xb_3$};}
    = (-1)^{\frac f2} \tilde{\chi}_{a,c}^f(\xb_1\xb_2^{-1},\xb_3).
    $$
The proper substitutions in this case are
    \begin{align*}
        x & \mapsto\mathrm{tr}(\xb_1\xb_2^{-1}) = t_1t_2-t_{12}, \\
        y & \mapsto\mathrm{tr}(\xb_3)=t_3, \\
        z & \mapsto\mathrm{tr}(\xb_3\xb_1\xb_2^{-1})=t_{13}t_2-t_{123}.
    \end{align*}

\subsection{Computations}
Recall that every rank 3 central function has 6 indices.  The first 4 come from the inclusion of a symmetric tensor $V_d \hookrightarrow V_a\otimes V_b \otimes V_c$.  Suppose it occurs with multiplicity $m$.  Then for $1\leq i,j\leq m$, the central functions are indexed by $\chi_{a,b,c,d}^{i,j}$, $i$ for $V_d$ and $j$ for $V^*_d$. Any 6-tuple arising in this fashion is called admissible. Note that this choice of index differs from $\chi_{a,b,c,d,e,f}$. By letting $c=0$ we recover the rank 2 central functions in \cite{LP} and by letting $b=c=0$ we recover the classical rank 1 central functions.

We call the set of all admissible 6-tuples satisfying $a+b+c=s$ the $s$-{\it order}. Table \ref{t:rank3strata} shows the first four orders of rank three central functions.

\begin{table}
0-order (1 function)

\begin{tabular}{>{$}l<{$}>{$}l<{$}}
 \chi_{0,0,0,0}^{1,1}=1 
\end{tabular}

\bigskip
1-order (3 functions)

\begin{tabular}{>{$}l<{$}>{$}l<{$}>{$}l<{$}}
\chi_{1,0,0,1}^{1,1}=t_1 & \chi_{0,1,0,1}^{1,1}=t_2 &\chi_{0,0,1,1}^{1,1}=t_3
\end{tabular}

\bigskip
2-order (9 functions)

\begin{tabular}{>{$}l<{$}>{$}l<{$}>{$}l<{$}}

\chi_{2,0,0,2}^{1,1}=t_1^2-1
    & \chi_{1,0,1,2}^{1,1}=\tfrac{t_1 t_3}{2}+\tfrac{t_{13}}{2}
        & \chi_{0,1,1,2}^{1,1}=\tfrac{t_2 t_3}{2}+\tfrac{t_{23}}{2}\\
\chi_{1,1,0,2}^{1,1}=\tfrac{t_1 t_2}{2}+\tfrac{t_{12}}{2}
    & \chi_{1,0,1,0}^{1,1}=t_1 t_3-t_{13}
        & \chi_{0,1,1,0}^{1,1}=t_2 t_3-t_{23}\\
\chi_{1,1,0,0}^{1,1}=t_1 t_2-t_{12} 
    & \chi_{0,2,0,2}^{1,1}=t_2^2-1
        & \chi_{0,0,2,2}^{1,1}=t_3^2-1
\end{tabular}

\bigskip
3-order (20 functions)

\begin{tabular}{>{$}l<{$}>{$}l<{$}}
\chi_{3,0,0,3}^{1,1}=t_1^3-2 t_1
    & \chi_{1,1,1,1}^{1,2}=-\frac{1}{2} t_1 t_2 t_3+\frac{t_{12} t_3}{2}+t_1 t_{23}-t_{123}\\
\chi_{2,1,0,3}^{1,1}=\frac{1}{3} t_2 t_1^2+\frac{2 t_{12}  t_1}{3}-\frac{2 t_2}{3}
    & \chi_{1,1,1,1}^{1,1}=\frac{3}{4} t_1 t_2 t_3+\frac{t_{12} t_3}{4}-\frac{t_2 t_{13}}{2}-\frac{t_1 t_{23}}{2}\\
\chi_{2,1,0,1}^{1,1}=t_2 t_1^2-t_{12} t_1-\frac{t_2}{2}
    & \chi_{1,0,2,3}^{1,1}=\frac{1}{3} t_1 t_3^2+\frac{2 t_{13} t_3}{3}-\frac{2 t_1}{3}\\
\chi_{2,0,1,3}^{1,1}=\frac{1}{3} t_3 t_1^2+\frac{2 t_{13} t_1}{3}-\frac{2 t_3}{3}
    & \chi_{1,0,2,1}^{1,1}=t_1 t_3^2-t_{13} t_3-\frac{t_1}{2}\\
\chi_{2,0,1,1}^{1,1}=t_3 t_1^2-t_{13} t_1-\frac{t_3}{2}
    & \chi_{0,3,0,3}^{1,1}=t_2^3-2 t_2\\
\chi_{1,2,0,3}^{1,1}=\frac{1}{3} t_1 t_2^2+\frac{2 t_{12} t_2}{3}-\frac{2 t_1}{3}
    & \chi_{0,2,1,3}^{1,1}=\frac{1}{3} t_3 t_2^2+\frac{2 t_{23} t_2}{3}-\frac{2 t_3}{3}\\
\chi_{1,2,0,1}^{1,1}=t_1 t_2^2-t_{12} t_2-\frac{t_1}{2}
    & \chi_{0,2,1,1}^{1,1}=t_3 t_2^2-t_{23} t_2-\frac{t_3}{2}\\
\chi_{1,1,1,3}^{1,1}=\frac{t_3 t_{12}}{3}+\frac{t_2 t_{13}}{3}+\frac{t_1 t_{23}}{3}
    & \chi_{0,1,2,3}^{1,1}=\frac{1}{3} t_2 t_3^2+\frac{2 t_{23} t_3}{3}-\frac{2 t_2}{3}\\
\chi_{1,1,1,1}^{2,2}=t_1 t_2 t_3-t_3 t_{12}
    & \chi_{0,1,2,1}^{1,1}=t_2 t_3^2-t_{23} t_3-\frac{t_2}{2}\\
\chi_{1,1,1,1}^{2,1}=\frac{1}{2} t_1 t_2 t_3-\frac{t_{12} t_3}{2}-t_2 t_{13}+t_{123}
    & \chi_{0,0,3,3}^{1,1}=t_3^3-2 t_3
\bigskip
\end{tabular}
\caption{Order of rank 3 central functions.}\label{t:rank3strata}
\end{table}

A more complicated example is
\begin{align*}
\chi_{3,2,2,3}^{2,1} &= \frac{1}{30} t_2^2 t_1^3+\frac{4}{15} t_2^2 t_3^2 t_1^3-\frac{1}{5} t_3^2 t_1^3-\frac{1}{5} t_2 t_3 t_{23}t_1^3+\frac{2 t_1^3}{15}-\frac{2}{15} t_2t_3^2 t_{12} t_1^2\\
                      &-\frac{7}{10} t_2^2 t_3t_{13} t_1^2+\frac{7}{15} t_3 t_{13}t_1^2+\frac{1}{5} t_3 t_{12} t_{23}t_1^2+\frac{3}{10} t_2 t_{13} t_{23}t_1^2+\frac{2}{3} t_2 t_3 t_{123}t_1^2\\
                      &-\frac{1}{3} t_{23} t_{123}t_1^2-\frac{1}{6} t_2^2 t_1+\frac{2}{15}t_2^2 t_3^2 t_1+\frac{1}{15} t_3^2t_1-\frac{2}{15} t_3^2 t_{12}^2t_1-\frac{1}{30} t_{12}^2t_1\\
                      &+\frac{2}{5}t_2^2 t_{13}^2 t_1-\frac{4}{15} t_{13}^2t_1+\frac{7}{30} t_{23}^2 t_1-\frac{1}{3} t_2t_3 t_{12} t_{13} t_1-\frac{1}{2} t_2 t_3t_{23} t_1\\
                      &+\frac{1}{30} t_{12} t_{13} t_{23}t_1+\frac{1}{3} t_3 t_{12} t_{123}t_1-\frac{1}{2} t_2 t_{13} t_{123}t_1-\frac{t_1}{15}+\frac{4}{15} t_2 t_{12}t_{13}^2\\
                      &+\frac{4}{15} t_2 t_3^2t_{12}-\frac{t_2 t_{12}}{10}+\frac{1}{30} t_3t_{12}^2 t_{13}+\frac{2}{15} t_2^2 t_3t_{13}-\frac{t_3 t_{13}}{10}-\frac{7}{30} t_3t_{12}t_{23}\\
                      &+\frac{2}{15} t_2 t_{13}t_{23}-\frac{1}{3} t_2 t_3t_{123}-\frac{1}{6} t_{12} t_{13}t_{123}+\frac{t_{23} t_{123}}{6}.
\end{align*}

Using both our tensorial algorithm and our combinatorial algorithm we have computed all orders up to 10 which gives 2254 known examples.  We note that on the same computer it took 2 minutes to compute all rank 3 central functions up to order 10 with the combinatorial algorithm, but took in excess of 6 hours to do so with the tensorial algorithm.

%% file: appendix.tex
\subsection{Some $\SL$ Representation Theory}\label{ss:reptheory}

We now review some basic $\SL$ representation theory, and reintroduce some notation from our first paper \cite{LP}.

Let $V_0=\C =V_0^*$ be the trivial representation of $\SL$.  Denote the standard basis for $\C^2$ by $e_1=\tmxt10$ and $e_2=\tmxt01$, and the dual basis by its transpose:
$e_1^*$ and $e_2^*.$

Then the standard representation and its dual are
$$V=V_1 = \C e_1\oplus \C e_2 \quad\text{and}\quad V^*=V_1^*=\C e_1^*\oplus \C e_2^*.$$
Denote the symmetric powers of these representations by
$$V_n=\mathrm{Sym}^n(V) \textrm{ and } V_n^*=\mathrm{Sym}^n(V^*).$$
One can show $V_n\approx (V_n)^*\approx V_n^*$ as $\SL$-modules, which is particular to $\SL$ and not obvious.

\begin{proposition}The symmetric powers of the standard representation of $\SL$ are all irreducible representations and
moreover they comprise a complete list.\end{proposition}

For proof see \cite{FH}.

The tensor product $V_a\otimes V_b,$ where $a,b \in \N $, is also a representation of $\SL$ and
decomposes into irreducible representations as follows:
\begin{proposition}[Clebsch-Gordan formula]\label{clebshgordan}
  \begin{equation}\label{eq:clebsch-gordan}
  V_a\otimes V_b\approx\bigoplus^{\mathrm{min}(a,b)}_{j=0} V_{a+b-2j}.
  \end{equation}
\end{proposition}

The particular $V_{a+b-2j}$ summands in this formula are described by the following:
\begin{definition}\label{def:admissible}
    Given $a,b\in\N$, we write $c\in\iadm{a,b}$ and say that $\{a,b,c\}$ is an \emph{admissible triple} for all $c=a+b-2j$, $0\le j\le \min(a,b)$.
\end{definition}

We remind ourselves of Schur's Lemma for later use:
\begin{proposition}[Schur's Lemma]
Let $G$ be a group, $V$ and $W$ irreducible representations of $G$, and $f\in\hm_G(V,W)$ with $f\neq0$.
If $V\approx W$, then $\dim_\C\mathrm{Hom}_G(V,W)=1$; and
if $V\not\approx W$, then $\dim_\C\mathrm{Hom}_G(V,W)=0$.
\end{proposition}

A tensor product $v_1\otimes v_2 \otimes \cdots \otimes v_n \in V^{\otimes n}$ projects to $V_n$ by \emph{symmetrizing}. We define its image under this operation by
    $$v_1\circ v_2 \circ \cdots \circ v_n \equiv \frac{1}{n!} \sum_{\sigma\in \Sigma_n} v_{\sigma(1)}\otimes v_{\sigma(2)} \otimes \cdots \otimes v_{\sigma(n)},$$
where the sum is over all permutations on $n$ elements. There exist bases for $V_n$ and $V_n^*$, given by
the elements%
\begin{align*}
    {\sf n}_{n-k}&=e_1^{n-k}e_2^k=\underbrace{e_1\circ e_1\circ \cdots \circ e_1}_{n-k}\circ
        \underbrace{e_2\circ e_2\circ \cdots \circ e_2}_{k}\quad\text{and}\\
    {\sf n}_{n-k}^*&=(e_1^*)^{n-k}(e_2^*)^k=\underbrace{e_1^*\circ e_1^*\circ \cdots \circ e_1^*}_{n-k}\circ
        \underbrace{e_2^*\circ e_2^*\circ \cdots \circ e_2^*}_{k}\:,
\end{align*}
respectively, where $0\le k\le n$. These elements are described in diagrammatic form in section \ref{ss:diagrammatic-symmetrization}.

The ``dual'' pairing between $V_n$ and $V_n^*$ is given by
$${\sf n}^*_{n-k}(v_1\circ v_2 \circ \cdots \circ v_n)=
\frac{1}{n!} \sum_{\sigma \in \Sigma_n}({\sf n}_{n-k})^*(v_{\sigma(1)}\otimes v_{\sigma(2)} \otimes
\cdots \otimes v_{\sigma(n)}),$$ where $\Sigma_n$ is the symmetric group on $n$ elements. In particular,
    $${\sf n}^*_{n-k}({\sf n}_{n-l})
    =\frac{(n-k)!k!}{n!}\delta_{kl}=\raisebox{2pt}{$\delta_{kl}$}\!\Big/\!\raisebox{-2pt}{$\tbinom{n}{k}$}.$$

Let $g=\tmx{g_{11} & g_{12}}{g_{21} & g_{22}}\in \SL$. The $\SL$-action on $V_n$ is given by
\begin{align*}
  g\cdot {\sf n}_{n-k}
    &=(g_{11}e_1+g_{21}e_2)^{n-k}(g_{12}e_1+g_{22}e_2)^k\\%
    &=\sum_{\substack{0\le j\le n-k\\0\le i \le k}}\tbinom{n-k}{j}\tbinom{k}{i}
        \left(g_{11}^{n-k-j}g_{12}^{k-i}g_{21}^{j}g_{22}^{i}\right){\sf n}_{n-(i+j)}.
\end{align*}

For the dual, $\SL$ acts on $V_n^*$ in the usual way:
\begin{displaymath}
    (g \cdot {\sf n}^*_{n-k})(v)= {\sf n}^*_{n-k}(g^{-1}\cdot v )
    \quad\textrm{ for } v\in V_n.
\end{displaymath} 

\subsection{$\SL$-character varieties}
Let $\C[x^k_{ij}]/\Delta$ be the complex polynomial ring in $4r$ variables ($1\leq k \leq r$ and $1\leq i,j \leq 2$), where $\Delta$ is the ideal generated by
the $r$ irreducible polynomials $$x^k_{11}x^k_{22}-x^k_{12}x^k_{21}-1.$$  It is not hard to see that $\C[\R_r]=\C[x^k_{ij}]/\Delta$.  Let $\widehat{x}^k_{ij}$ be the image of $x^k_{ij}$ under $\C[x^k_{ij}] \to \C[x^k_{ij}]/\Delta$.
Define $$\xb_k=\left(
\begin{array}{cc}
\widehat{x}^k_{11} & \widehat{x}^k_{12} \\
\widehat{x}^k_{21} & \widehat{x}^k_{22}  \\
\end{array}\right).$$  Such elements are called {\it generic unimodular matrices}.  We note that $$\left(
\begin{array}{cc}
x^k_{11} & x^k_{12} \\
x^k_{21} & x^k_{22}  \\
\end{array}\right)$$ are simply called {\it generic matrices}.
Let $G=\SL$.  $G$ acts on $\C[\R_r]$ as follows:  $$g\cdot \widehat{x}^k_{ij}=y^k_{ij} \text{ where } \yb_k:=g\xb_k g^{-1},\text{ and }g\in G.$$
The ring of invariants $\C[\R_r]^{G}$ is a finitely-generated domain (see \cite{Na}), which implies $$\X_r:=\mathrm{Spec}_{max}\left(\C[\R_r]^{G}\right)$$ is an affine variety over $\C$, called the
$G$-{\it character variety} of $\F_r$.  It is the variety whose coordinate ring is the ring of invariants, that is, $$\C[\X_r]=\C[\R_r]^{G}.$$

Closely related to $\C[\X_r]$ is the ring of invariants $$\C[\mathfrak{Y}_r]:=\C[\mathfrak{gl}(2,\C)^{\times r}]^{\SL}=\C[x^k_{ij}]^{\SL}.$$
In fact one can show (see \cite{L3})
$$\C[\mathfrak{Y}_r]/\Delta\approx \C[\X_r]$$
Otherwise stated, $$\C[x^k_{ij}]^{\SL}/\Delta\approx \left( \C[x^k_{ij}]/\Delta \right)^{\SL},$$ which is true because $\SL$ is {\it linearly} reductive.
In 1976 Procesi proved (in the context of $n\times n$ generic matrices)
\begin{theorem}[Procesi]
$\C[\mathfrak{Y}_r]$ is generated by the invariants $$\{\Tr{\xb_{i_1}\xb_{i_2}\cdots\xb_{i_k}}\},$$ where $\xb_j$ are generic matrices.
\end{theorem}
Evidently, this ring is multigraded.
Finding minimal generators amounts to finding all linear relations among generators of the same multidegree in the vector space
$$\mathcal{V}_r=\C[\mathfrak{Y}_r]^+/\left(\C[\mathfrak{Y}_r]^+\right)^2,$$ where $\C[\mathfrak{Y}]^+$ is the ideal of positive terms.
It is worth observing that there can be no relation among generators of differing multidegree.  Any such linear relation in $\mathcal{V}_r$ pulls back to a reduction relation in $\C[\Y_r]$ which in turns projects to a ``reduction relation'' in $\C[\X_r]$.  The Cayley-Hamilton equation gives $$\xb^2-\Tr{\xb}\xb+\mathrm{det}(\xb)\id=0.$$  And if we assume $\mathrm{det}(\xb)=1$, as is the case in $\C[\X_r]$, we easily derive $\Tr{\xb^{-1}}=\Tr{\xb}$ and $\Tr{\xb^2}=\Tr{\xb}^2-2$.  Hence the generators $\Tr{\xb^2}$ in $\C[\Y_r]$ project to $\Tr{\xb}^2-2$ in $\C[\X_r]$ and so are freely eliminated.

We work a couple examples before moving on:
$\C[\Y_1]$ is algebraically generated by $\{\Tr{\xb^n}\}$.  But if $n\geq 3$ the  Cayley-Hamilton equation provides relations which
express the generator in terms of the generators $\Tr{\xb}$ and $\Tr{\xb^2}$.  However, the dimension of this variety is computed to be $2$.  Thus there can be no further relations.

Let $\mathcal{M}(R)$ be a set of minimal generators of a ring $R$.  Then in \cite{L3} it is shown, as with a generalization to $\SLm{n}$, that $$\mathcal{M}(\C[\Y_r])-\{\Tr{\xb_1^2},...,\Tr{\xb_r^2}\}=\mathcal{M}(\C[\X_r]),$$ as long as the generators of $\C[\Y_r]$ are taken to be of the form in Procesi's theorem.

Multiplying the Cayley-Hamilton equation on both sides by words $\ub$ and $\vb$ allows us to freely eliminate the generators of type: $\Tr{\ub\xb^n\vb}$ as long as $n\geq 2$ and at least one of $\ub$ or $\vb$ is not the identity.
So for the case $\C[\Y_2]$, we are left with the generators $$\{\Tr{\xb_1},\Tr{\xb_2},\Tr{\xb_1^2},\Tr{\xb_2^2},\Tr{\xb_1\xb_2}\}$$ since
any other expression in two letters would result in a sub-expression with an exponent greater than one, which we just showed was impossible.
Since in this case the dimension of the variety is $5$, there can be no further relations and thus these generators are minimal and $\Y_2\approx \C^5$.  We can conclude that $\X_1\approx \C$ and $\X_2\approx \C^3$.

More generally, it can be shown that there are no generators necessary that have word length 4 or more (see \cite{G9} for an exposition).  In particular,
\begin{align}
\Tr{\wb_1\wb_2\wb_3\wb_4}=&\Tr{\wb_1}\Tr{\wb_2}\Tr{\wb_3}\Tr{\wb_4}+\Tr{\wb_1}\Tr{\wb_2\wb_3\wb_4}\nonumber\\
&+\Tr{\wb_2}\Tr{\wb_3\wb_4\wb_1}+\Tr{\wb_3}\Tr{\wb_4\wb_1\wb_2}\nonumber\\
&+\Tr{\wb_4}\Tr{\wb_1\wb_2\wb_3}+\Tr{\wb_1\wb_2}\Tr{\wb_3\wb_4}\nonumber\\
&+\Tr{\wb_1\wb_4}\Tr{\wb_2\wb_3}+\Tr{\wb_1\wb_2}\Tr{\wb_3\wb_4}\nonumber\\
&-\Tr{\wb_1\wb_3}\Tr{\wb_2\wb_4}-\Tr{\wb_1}\Tr{\wb_2}\Tr{\wb_3\wb_4}\nonumber\\
&-\Tr{\wb_1}\Tr{\wb_4}\Tr{\wb_2\wb_3}-\Tr{\wb_2}\Tr{\wb_3}\Tr{\wb_1\wb_4}\nonumber\\
&-\Tr{\wb_3}\Tr{\wb_4}\Tr{\wb_1\wb_2}\label{length4reduction}.
\end{align}
and
\begin{align}
\Tr{\xb\yb\zb}+\Tr{\xb\zb\yb}&=\Tr{\xb\yb}\Tr{\zb}+\Tr{\xb\zb}\Tr{\yb}\nonumber\\
&+\Tr{\zb\yb}\Tr{\xb}-\Tr{\xb}\Tr{\yb}\Tr{\zb}.\label{sumformula}\end{align}

Equations \eqref{length4reduction} and \eqref{sumformula} together imply there are $N_r=\frac{r(r^2+5)}{6}$ minimal generators of $\C[\X_r]$.  They are:
$\mathcal{G}_1=\{\Tr{\xb_1},...,\Tr{\xb_r}\}$ of order $r$,
$\mathcal{G}_2=\{\Tr{\xb_i\xb_j}\ |\ 1\leq i<j \leq r\}$ of order $\frac{r(r-1)}{2}$, and
$\mathcal{G}_3=\{\Tr{\xb_i\xb_j\xb_k}\ | \ 1\leq i<j<k \leq r \}$ of order $\frac{r(r-1)(r-2)}{3}$.

Enumerating these $N_r$ generators, $\{t_1,...,t_{N_r}\}$ defines a polynomial mapping
$$t=(t_1,...,t_{N_r}):\X_r\longrightarrow \C^{N_r}.$$  It is not hard to show that $t$ is a proper injection and hence defines a homeomorphism (with respect to the induced topology from $\C^{N_r}$) onto its image.  We conclude with the following global geometric result:  the smallest affine embedding $\X_r\longrightarrow \C^{N}$ is when $N=\frac{r(r^2+5)}{6}$, which follows from the fact that $\{t_1,...,t_{N_r}\}$ is a {\it minimal} generating set.

In Section \ref{s:rankthree} we discussed that $\X_3$ is a branched double cover of $\C^6$.  In general, $$\X_r=\mathrm{Spec}_{max}\left(\C[t_1,...,t_{N_r}]/\mathfrak{I}_r\right)$$ where $\mathfrak{I}_r$ is an ideal generated by $\frac{1}{2}\left(\binom{r}{3}^2+\binom{r}{3}\right)+r\binom{r}{4}$ polynomials (see \cite{Dr}).